

\magnification=1200
\pageno=1
\baselineskip=18pt
\parindent=0pt

\font\bigbf=cmbx10 scaled\magstep1
\raggedbottom

\input eplain

\newfam\bbbfam
\font\bbbten=msbm10
\font\bbbseven=msbm7
\font\bbbfive=msbm5
\textfont\bbbfam=\bbbten
\scriptfont\bbbfam=\bbbseven
\scriptscriptfont\bbbfam=\bbbfive
\def\bbb{\fam=\bbbfam}

\font\teneufm=eufm10
\font\seveneufm=eufm7
\font\fiveeufm=eufm5
\newfam\eufmfam
\textfont\eufmfam=\teneufm
\scriptfont\eufmfam=\seveneufm
\scriptscriptfont\eufmfam=\fiveeufm
\def\eufm#1{{\fam\eufmfam\relax#1}}

\font\teneufb=eufb10
\font\seveneufb=eufb7
\font\fiveeufb=eufb5
\newfam\eufbfam
\textfont\eufbfam=\teneufb
\scriptfont\eufbfam=\seveneufb
\scriptscriptfont\eufbfam=\fiveeufb

\font\teneurm=eurm10
\font\seveneurm=eurm7
\font\fiveeurm=eurm5
\newfam\eurmfam
\textfont\eurmfam=\teneurm
\scriptfont\eurmfam=\seveneurm
\scriptscriptfont\eurmfam=\fiveeurm

\font\teneurb=eurb10
\font\seveneurb=eurb7
\font\fiveeurb=eurb5
\newfam\eurbfam
\textfont\eurbfam=\teneurb
\scriptfont\eurbfam=\seveneurb
\scriptscriptfont\eurbfam=\fiveeurb

\rm
\centerline{\bigbf Triple product integrals and Rankin-Selberg {\it L\/}-functions}
 \vskip10pt
\centerline{by Andr\'as BIR\'O\footnote{}{Research partially
supported by NKFIH (National Research, Development and Innovation Office)
grants K135885, K119528, K143876 and by the R\'enyi Int\'ezet Lend\"ulet Automorphic Research Group}}
 \footnote{}{} \footnote{}{2020
}\hfill\break
\centerline{A. R\'enyi Institute of Mathematics}

\centerline{ 1053 Budapest, Re\'altanoda u. 13-15., Hungary; e-mail: biro.andras@renyi.hu}
 \vskip20pt

\noindent {\bf Abstract}. We prove a reciprocity formula that relates a
spectral average of products of triple product integrals involving automorphic forms of
weights $0$ and $1/2$ to the classical Rankin-Selberg integrals for automorphic
forms of weight $0$.

\vskip20pt

\noindent{\bf 1. Introduction }
\medskip

{\bf 1.1. Triple product integrals of weight} $0$ {\bf and weight} $ ${\bf 1/2 }
{\bf Maass forms.} Let $u$ and $U$ be two Maass cusp forms of weight $
0$ for
$SL(2,{\bf Z})$. This means that $u$ and $U$ are $SL(2,{\bf Z})$-invariant
functions on the open upper half plane $\bbb H$ \xrdef{H}
decaying exponentially as $\hbox{\rm Im$\,z\rightarrow\infty$}$, and $
u$ and $U$ are eigenfunctions
of the hyperbolic Laplace operator $\Delta_0:=y^2\left({{\partial^
2}\over {\partial x^2}}+{{\partial^2}\over {\partial y^2}}\right)$. Let
$D_1$ \xrdef{D1} be a fundamental domain of the quotient
$SL(2,{\bf Z})\setminus\bbb H$ and $d\mu_z:={{dxdy}\over {y^2}}$. \xrdef{dmuz} We write
$$\left(f_1,f_2\right)_1:=\int_{D_1}f_1(z)\overline {f_2(z)}d\mu_z
,$$
 \xrdef{f1f2}
where $d\mu_z$ is an $SL(2,{\bf R})$-invariant measure on $\bbb H$. The triple product integral
$$\left(\left|U\right|^2,u\right)_1\eqno (1.1)$$
is an important object of study in the theory of automorphic
forms. For example, the famous Quantum Unique
Ergodicity (QUE) Conjecture states that if $u$ is fixed, $U$ is a Hecke
eigenform satisfying $\left(U,U\right)_1=1$  and the Laplace eigenvalue of
$U$ tends to $-\infty$, then (1.1) tends to $0$. This conjecture was proved by Lindenstrauss and Soundararajan (see [L]
and [S]). However, quantifying the rate of
convergence in QUE is still an open problem. Watson (see
[Wa]) proved an important identity relating (1.1) to the
central value of a degree 8 $L$-function. This
identity shows that the Generalized Riemann Hypothesis for some $GL(2)\times GL(3)$ Rankin-Selberg $L$-functions would give a quantitative form of QUE.

The integrals (1.1) can be expressed in terms of triple
product integrals involving weight $1/2$ Maass forms, see [B1], Theorem 1.1. This
motivates the study of the triple product integrals to be
considered in Theorem 1.1 below. To define them properly
and to state our main result we need some notations.

{\bf 1.2. Necessary notations.} There is a
list of notations at the end of the paper.

We write
$$\Gamma_0(4):=\left\{\left(\matrix{a&b\cr
c&d\cr}
\right)\in SL(2,{\bf Z}):\hbox{\rm \ $c\equiv 0\hbox{\rm \ (mod $
4$})$}\right\}.$$
 \xrdef{Gammazero4}
Let $D_4$ \xrdef{D4} be a fundamental domain of the quotient
$\Gamma_0(4)\setminus\bbb H$ and
$$\left(f_1,f_2\right)_4:=\int_{D_4}f_1(z)\overline {f_2(z)}d\mu_z
.$$
 \xrdef{f1f24}
The hyperbolic Laplace operator of weight $l$ is given by:
$$\Delta_l:=y^2\left({{\partial^2}\over {\partial x^2}}+{{\partial^
2}\over {\partial y^2}}\right)-ily{{\partial}\over {\partial x}}.$$
 \xrdef{Laplace}
For a complex number $z\neq 0$, its argument (denoted by
$\arg z$) \xrdef{arg} is chosen to be in the range $(-\pi ,\pi ]$, and we define
$\log z:=\log\left|z\right|+i\arg z$ and $z^s:=e^{s\log z}$ for any $
s\in {\bf C}$.

We write $e(x):=e^{2\pi ix}$. \xrdef{e(x)} For $z\in\bbb H$, we define
$$B_0(z):=\left(\hbox{\rm Im$\,z$}\right)^{{1\over 4}}\theta\left(z\right
)=\left(\hbox{\rm Im$\,z$}\right)^{{1\over 4}}\sum_{m=-\infty}^{\infty}
e(m^2z).\eqno (1.2)$$
 \xrdef{theta}
 \xrdef{B0}
We define the symbol $\left({c\over d}\right)$ where $c$ is an integer and $d$ is an odd integer. For $d>0$ this is the usual Jacobi symbol, and
we extend it by the formulas $\left({c\over d}\right):={c\over {\left
|c\right|}}\left({c\over {-d}}\right)$ for $c\neq 0$,
$\left({0\over d}\right):=1$ for $d=\pm 1$, $\left({0\over d}\right
):=0$ for $\left|d\right|>1$. Define $\epsilon_d:=1$ for
$d\equiv 1$ $\left(\hbox{\rm mod 4}\right)$, $\epsilon_d:=i$ for $
d\equiv -1$ $\left(\hbox{\rm mod 4}\right)$. For
$\gamma =$$\left(\matrix{a&b\cr
c&d\cr}
\right)\in\Gamma_0(4)$ let $\nu (\gamma ):=\left({c\over d}\right)\overline {
\epsilon_d}.$ \xrdef{nu}

Then for every $z\in\bbb H$ and $\gamma\in\Gamma_0(4)$ we have
$$B_0(\gamma z)=\nu (\gamma )\left({{j_{\gamma}(z)}\over {\left|j_{
\gamma}(z)\right|}}\right)^{1/2}B_0(z),\eqno (1.3)$$
where $j_{\gamma}(z):=cz+d$ for $\gamma =$$\left(\matrix{a&b\cr
c&d\cr}
\right)\in SL(2,{\bf R})$. \xrdef{jgamma} It
is also known that $B_0$ has an additional
transformation formula
$$B_0\left({{-1}\over {4z}}\right)=e\left({{-1}\over 8}\right)\left
({z\over {\left|z\right|}}\right)^{{1\over 2}}B_0(z)\eqno (1.4)$$
for every $z\in\bbb H$.

In this paper, any automorphic function is of weight
$l={1\over 2}+2n$ or $l=2n$ with some integer $n$. A smooth function $
f:\bbb H\rightarrow {\bf C}$ is said to be an automorphic
function of weight $l$ for $\Gamma$ if it has at most polynomial
growth at the cusps of $\Gamma$ and satisfies the transformation formula
$$f(\gamma z)=\left({{j_{\gamma}(z)}\over {\left|j_{\gamma}(z)\right
|}}\right)^lf(z)\cdot\left\{\matrix{1&\quad {\rm i}{\rm f}&l=2n\cr
\cr
\nu (\gamma )&\quad {\rm i}{\rm f}&\qquad l=2n+{1\over 2}\cr}
\right.$$
for any $z\in\bbb H$ and $\gamma\in\Gamma$, where $\Gamma$ is either $
SL(2,{\bf Z})$ or
$\Gamma_0(4)$. If $l={1\over 2}+2n$, we can take only
$\Gamma =\Gamma_0(4)$. The operator $\Delta_l$ acts on automorphic
functions of weight $l$. We say that $f$ is a Maass form of
weight $l$ for $\Gamma$, if $f$ is an automorphic function and it is an eigenfunction of $
\Delta_l$. If a Maass form $f$ has exponential decay at
all of the cusps of $\Gamma$, it is called a cusp form.

A Maass cusp form $f$ of weight $0$ for $SL(2,{\bf Z})$ is called
{\it even\/} if $f\left(z\right)=f\left(-\overline z\right)$, and it is called {\it odd\/} if
$f\left(z\right)=-f\left(-\overline z\right)$.

In this work, our weight $0$ Maass cusp forms $u_1$, $u_2$ \xrdef{u1u2}
for $SL(2,{\bf Z})$ are assumed to be

(i) $L^2-$normalized (i.e. $\left(u_j,u_j\right)_1=1$ for $j=1,2)$,

(ii) either orthogonal to each other (i.e. $\left(u_1,u_2\right)_
1=0$) or satisfying $u_1=u_2$,

(iii) and either both even or both odd.

Assume that $\Delta_0u_j=s_j(s_j-1)u_j$, where $s_j={1\over 2}+it$$_
j$ and $t_j>0$ \xrdef{tj}
($j=1,2$). We have the Fourier expansions \xrdef{rhou1rhou2}
$$u_1(z)=\sum_{m\neq 0}\rho_{u_1}(m)W_{0,it_1}(4\pi\left|m\right|
y)e(mx),\;u_2(z)=\sum_{m\neq 0}\rho_{u_2}(m)W_{0,it_2}(4\pi\left|
m\right|y)e(mx).\eqno (1.5)$$
Here $W_{\alpha ,\beta}$ denotes the Whittaker functions, see
Section 3.7 for the definition of these functions.

The Rankin-Selberg $L$-function is defined in terms of an
absolutely convegent Dirichlet series
$$L\left(S\right)=L\left(S,u_1\otimes\overline {u_2}\right):=\sum_{
m>0}\rho_{u_1}(m)\overline {\rho_{u_2}(m)}m^{1-S}\eqno (1.6)$$
 \xrdef{Rankin}
for $\hbox{\rm Re$\, S$}\gg 1$. It is well-known that $L\left(S\right)$ extends meromorphically to the whole complex
plane and is regular for $\hbox{\rm Re$\, S$}\ge 1/2$ with at most a simple
pole at $S=1$. Such a simple pole occurs only when
$u_1=u_2$.

The Wilson function $\phi_{\lambda}\left(x;a,b,c,d\right)$ was defined in [G1],
we give its definition in Section 3.6. We use the
abbreviations  $\Gamma\left(X\pm Y\right):=\Gamma\left(X+Y\right)\Gamma\left
(X-Y\right)$ and
$$\Gamma\left(X\pm Y\pm Z\right):=\Gamma\left(X+Y+Z\right)\Gamma\left
(X+Y-Z\right)\Gamma\left(X-Y+Z\right)\Gamma\left(X-Y-Z\right).$$
 \xrdef{gamma}
Recalling the notations $t_j$ and $s_j$ ($j=1,2$) from above define
$$N^{+}\left(S,t\right):={{\Gamma\left(S\pm it_1+it_2\right)\Gamma\left
({1\over 4}+it_2\pm it\right)}\over {\sin\pi\left(2it_2\right)}}\left
(\sin\pi s_1+\sin\pi\left(1-s_2-S\right)\right)\phi_{i\left
({1\over 2}-S\right)}^{+}\left(t\right),$$
$$N^{-}\left(S,t\right):={{\Gamma\left(S\pm it_1-it_2\right)\Gamma\left
({1\over 4}-it_2\pm it\right)}\over {\sin\pi\left(-2it_2\right)}}\left
(\sin\pi s_1+\sin\pi\left(s_2-S\right)\right)\phi_{i\left
({1\over 2}-S\right)}^{-}\left(t\right),$$
$$\phi_{\lambda}^{+}\left(x\right):=\phi_{\lambda}\left(x;{3\over
4}+it_2,{1\over 4}+it_1,{1\over 4}-it_1,{3\over 4}-it_2\right),$$
$$\phi_{\lambda}^{-}\left(x\right):=\phi_{\lambda}\left(x;{3\over
4}-it_2,{1\over 4}+it_1,{1\over 4}-it_1,{3\over 4}+it_2\right),$$
and let \xrdef{phi+}
$$N\left(S,t\right):=N^{+}\left(S,t\right)+N^{-}\left(S,t\right).$$
 \xrdef{N(s,t)}
This function was introduced in [B2].

{\bf CONVENTION.} Since the Maass cusp forms $u_1$, $u_2$ and the positive numbers $
t_1$ and $t_2$ are
fixed, we will not denote the dependence on $t_1$ and $t_2$ in the
sequel.

{\bf 1.3. The main result.}

Denote by $L_l^2(D_4)$ \xrdef{Ll(d4)} the space of automorphic functions of weight $
l$$ $ for $\Gamma_0(4)$ for
which $\left(f,f\right)_4<\infty$.

Take $u_{0,1/2}=c_0B_0$, where $c_0$ is chosen such that
$\left(u_{0,1/2},u_{0,1/2}\right)_4=1$. Let $\left\{u_{j,1/2}:\;j
\ge 0\right\}$ \xrdef{u(j,1/2)} be an orthonormal basis of Maass
forms for the discrete part of $L_{1/2}^2$$(D_4)$. Write
$$\Delta_{1/2}u_{j,1/2}=\Lambda_ju_{j,1/2},\;\Lambda_j=S_j(S_j-1)
,\;S_j={1\over 2}+iT_j.$$
 \xrdef{Tj}
It is known that $\Lambda_0=-{3\over {16}}$ and $\Lambda_j\rightarrow
-\infty$. It follows from
[Sa], Theorem 3.6 that $\Lambda_j<-{3\over {16}}$ for $j\ge 1$.

For the cusps $\eufm{a}=0,\infty$ denote by $E_{\eufm{a}}\left(z,
s,{1\over 2}\right)$
the Eisenstein series of weight ${1\over 2}$ for the group $\Gamma_
0(4)$ at the
cusp $\eufm{a}$. We give its definition for $z\in\bbb H$ and
Re$\, s>1$ in Section 2.5. On the one hand, as a function of $z$ it is an eigenfunction of $
\Delta_{1/2}$ of eigenvalue
$s(s-1)$. On the other hand, for every $z$ the function $E_{\eufm{
a}}\left(z,s,{1\over 2}\right)$ has a meromorphic continuation
in $s$ to the whole plane, and this function is regular at
every point $s$ with Re$\, s={1\over 2}$. If $f$ is an automorphic function of weight $
1/2$ and the following
integral is absolutely convergent, define
$$\zeta_{\eufm{a}}(f,r):=\int_{D_4}f(z)\overline {E_{\eufm{a}}\left
(z,{1\over 2}+ir,{1\over 2}\right)}d\mu_z.$$
 \xrdef {zeta(a,f,r)}
Let $\beta >0$. We say that a function $\chi$ satisfies condition $
C_{\beta}$ \xrdef {Cbeta}  if $\chi$ is an even holomorphic function defined on the
strip $\left|\hbox{\rm Im}\,z\right|<\beta$ and for every fixed $K>
0$ the function
$$\left|\chi (z)\right|e^{-\pi\left|z\right|}\left(1+\left|z\right
|\right)^K$$
is bounded on this strip.

Let $\delta_{u_1,u_2}$ be Kronecker's \xrdef {kron} symbol. We write
$\left(\kappa (u)\right)(z):=u(4z)$. \xrdef {kappa} We denote by $\zeta\left(S\right
)$ the Riemann zeta
function. \xrdef {zeta}

{\bf THEOREM 1.1.} {\it There is an absolute constant} $\beta >0$ {\it such that if} $
\chi$ {\it is a function
satisfying condition} $C_{ \beta}${\it , then the sum of}
$$\sum_{j=1}^{\infty}\chi\left(T_j\right)\left(B_0\kappa\left(\overline {
u_2}\right),u_{j,{1\over 2}}\right)_4\overline {\left(B_0\kappa\left
(\overline {u_1}\right),u_{j,{1\over 2}}\right)_4}\eqno (1.7)$$
 {\it
and}
$${1\over {4\pi}}\sum_{\eufm{a}=0,\infty}\int_{-\infty}^{\infty}\chi\left
(r\right)\zeta_{\eufm{a}}\left(B_0\kappa\left(\overline {u_2}\right
),r\right)\overline {\zeta_{\eufm{a}}\left(B_0\kappa\left(\overline {
u_1}\right),r\right)}dr\eqno (1.8)$$
 {\it equals the sum of }
$${3\over {2\pi^{3/2}}}{{\delta_{u_1,u_2}}\over {\Gamma\left({1\over
2}\pm it_1\right)}}\int_{-\infty}^{\infty}{{\Gamma\left({1\over 4}
\pm it\right)\Gamma\left({1\over 4}\pm it\pm it_1\right)}\over {\Gamma\left
(\pm 2it\right)}}\chi (t)dt$$
 {\it and}
$$-{6\over {\Gamma\left({1\over 2}\pm it_2\right)}}{1\over {2\pi
i}}\int_{\left({1\over 2}\right)}{{\zeta\left(2S\right)L\left(S\right
)}\over {\left(2\pi\right)^{2S}}}\Gamma\left(S\right)\Gamma\left(
1-S\right)H_{\chi}\left(S\right)dS,\eqno (1.9)$$
 {\it where}
$$H_{\chi}\left(S\right):=\int_{-\infty}^{\infty}{{\Gamma\left({1\over
4}\pm it\right)\Gamma\left({1\over 4}\pm it\pm it_1\right)}\over {
\Gamma\left(\pm 2it\right)}}\chi (t)N\left(S,t\right)dt.$$
 \xrdef{Hchi(S)}

{\it The sum in (1.7), and the integrals in (1.8) and (1.9) are
absolutely convergent.}

{\bf 1.4. Discussion of the main result.}

{\bf REMARK 1.1.} Many ideas of our proof are present also in
papers of Nelson, see [Nel1], [Nel2]. See, in
particular, [Nel2, formula (10)], the discussion below that
formula and [Nel2, formulas (14), (11)]. Indeed, using our notation, Nelson considered the following quantities:
$$\left|\int_{D_4}B_0\left(z\right)\phi\left(z\right)\overline {h\left
(z\right)}d\mu_z\right|^2,\eqno (1.10)$$
where $\phi\left(z\right)$ and $h\left(z\right)$ are cusp forms for $
\Gamma_0(4)$ of
weights 0 and $1/2$, respectively. He suggested summing (1.10) over either $\phi$ or $h$ in an orthonormal basis, and then expressing the resulting sum using Parseval's identity as an inner product involving $|B_{0}|^2$, i.e.,
$$\int_{D_4}\left|B_0\left(z\right)\right|^2\left|\phi\left(z\right
)\right|^2d\mu_z.\eqno (1.11)$$
Then he remarks in [Nel2, formula (14)] that $\left|B_0\right|^2$ is orthogonal to cusp forms, which implies that $\left
|B_0\right|^2$ can
be expressed as a linear combination of Eisenstein
series, see [Nel2, formula (11)]. Then one can unfold the
integral (1.11), and this leads to Rankin-Selberg $L$-functions.

In this paper, rather than simply summing (1.10) over $h$, we insert a weight function that depends on its Laplace eigenvalue. Although Parseval's identity cannot be applied in this case, the resulting sum can still be expressed as a sum of inner products involving $B_{0} \overline{B_{n}}$, where the functions $B_{n}$ are liftings of $B_{0}$ via the Maass operators. These products are still linear combinations
of Eisensein series, a technical variant of this key fact is proved in Lemma 4.7
below. Then we can apply the unfolding method,
getting again an expression involving Rankin-Selberg
$L$-functions.

Many convergence problems occur during this process,
but finally we are able to give an explicit class of
admissible test functions and an explicit form of
the integral transform.

Also, instead of the absolute square in (1.10) we consider
the product of two such triple product integrals with two different weight $0$ Maass cusp forms $u_1$, $u_2$ for $SL(2,{\bf Z})$.

Note that the fact that $\left|B_0\right|^2$ is orthogonal to cusp
forms played a role already in our work [B3] (see Lemma 6.6
there), where a duality relation was proved for the kind of inner products
considered also in this paper. The duality relation proved
in [B3] involved also holomorphic analogues of the triple
product integrals of Theorem 1.1 above. It is possible to
prove an analogue of Theorem 1.1 also for such inner
products. We will state this holomorphic analogue
without proof in Section 1.5.

We will give a bit more detailed sketch of the proof of
Theorem 1.1 in Section 1.6.

{\bf REMARK 1.2.} In this remark we show that it is
reasonable to expect that a special case of our formula
recovers a particular instance of the spectral reciprocity formulae discovered recently by Humphries-Khan and Kwan in [H-K] and [Kw].

Remark 1.2 can be skipped, the rest of the paper can be understood without reading it. Some notions involved in the present remark will not be used later in the paper, therefore instead of giving every
definition here we just refer to the literature. Our main references will be [B1] and [K-S], most of the notions are defined there.

We will consider the cuspidal sum (1.7) of our Theorem 1.1 above
in the case $u_1=u_2$,  and assume also that $u_1$ is a simultaneous Hecke
eigenform. We first choose our orthonormal basis
$\left\{u_{j,1/2}:\;j\ge 0\right\}$ in a special way. In order to do that we
have to define some operators.

The Hecke operator $T_{p^2}$ \xrdef{Tp2} of weight ${1\over 2}$ for every prime
$p\neq 2$ and the operator $L$ \xrdef{l} are defined in [K-S], pp 199-200
and p 195, respectively. These operators act on the space $L_{1/2}^
2$$(D_4)$,
they are self-adjoint and commute with each other and with $\Delta_{
1/2}$.
Hence our orthonormal basis
$\left\{u_{j,1/2}:\;j\ge 0\right\}$ can be chosen in such a way that every
$u_{j,1/2}$ is an eigenfunction of the operators $T_{p^2}$ ($p\neq 2$) and of the operator $L$ (see [K-S],
pp 195-196). By Lemma 5.3 and Lemma 5.5 (ii) of [B1] we
see that $B_0\kappa\left(\overline {u_1}\right)$ is an eigenfunction of $
L$ of eigenvalue $1$.
But two $L$-eigenfunctions with different $L$-eigenvalues
are orthogonal to each other. Therefore we can keep in
(1.7) only those $u_{j,1/2}$ having $L$-eigevalue 1, since the
contribution of other terms is 0. For the case
$Lu_{j,1/2}=u_{j,1/2}$ we will prove Proposition 1.1 below. We
first need some notations.

Let $F\in\hbox{\rm $L_{1/2}^2$$(D_4)$}$ be a cusp form of weight $
1/2$ for $\Gamma_0(4)$
which is an eigenfunction of the Hecke operator $T_{p^2}$ of
weight ${1\over 2}$ for every prime $p\neq 2$ and satisfies $LF=F$.
Assume also $\left(F,F\right)_4=1$. Assume $\rho_F(1)\neq 0$, where $
\rho_F(1)$ is the first Fourier
coefficient of $F$ at $\infty$. Under this assumption the Shimura lift $\hbox{\rm Shim$
F$}$ is defined \xrdef {Shim}
in [K-S], pp. 196-197. It is an even Maass cusp form of weight $0$
for ${\rm S}{\rm L}(2,{\rm Z})$, it is a simultaneous Hecke eigenform and its
first Fourier coefficient is $1$. Let $U$ be a cusp form and a simultaneous
Hecke eigenform of weight $0$ for ${\rm S}{\rm L}(2,{\rm Z})$ satisfying $\left
(U,U\right)_1=1$.

Assume $\Delta_0\left(\hbox{\rm Shim$F$}\right)=\left(-{1\over 4}
-t^2\right)\hbox{\rm Shim$F$},$ $\Delta_0U=\left(-{1\over 4}-T^2\right
)U,$
$\Delta_{1/2}F=\left(-{1\over 4}-r^2\right)F$. Note that we have $
t=2r$ e.g. by Theorem 1 of [B4].

{\bf PROPOSITION 1.1.} {\it Assume that} $\rho_F(1)\neq 0$. {\it Using the notations and assumptions above we have that} $\left
|\left(B_0\kappa\left(\overline U\right),F\right)_4\right|^2$ {\it equals}
$$d{{\left|\rho_U\left(1\right)\right|^2L\left({1\over 2},\hbox{\rm Shim$
F$}\otimes\hbox{\rm $\hbox{\rm sym}^2U$}\right)}\over {\left(\hbox{\rm Shim$
F$},\hbox{\rm Shim$F$}\right)_1}}\left|\Gamma\left({{{1\over 2}+i
t}\over 2}\right)\right|^2\left|\Gamma\left({{{1\over 2}+2iT\pm i
t}\over 2}\right)\right|^2,\eqno $$
{\it where} $\hbox{\rm sym}^2U$ {\it is the symmetric square lift of} $
U$, $L\left(s,\hbox{\rm Shim$F$}\otimes\hbox{\rm $\hbox{\rm sym}^
2U$}\right)$
{\it is the Rankin-Selberg} $L${\it -function of the pair} $\left
(\hbox{\rm Shim$F$, }\hbox{\rm sym}^2U\right)${\it , and} $d>0$ {\it is an absolute constant}.

The Shimura lift $\hbox{\rm Shim$F$}$ is defined also without the
condition $\rho_F(1)\neq 0$ on p 981 of [D-I-T]. It is very likely
that using that definition Proposition 1.1 is true without the
condition $\rho_F(1)\neq 0$, but we were able to prove it only
under this condition.

Assume now that Proposition 1.1 is true without the
condition $\rho_F(1)\neq 0$. Let $u_1=u_2$,  and assume also that $
u_1$ is a simultaneous Hecke
eigenform. We can then see that choosing the test
functions suitably the cuspidal sum (1.7) of Theorem 1.1
above coincides with the cuspidal sum of
Theorem 1.1 of [Kw] assuming there that $s={1\over 2}$ and $\Phi$ is self-dual.

Indeed, we choose $U=u_1$ in Proposition 1.1. Then $U$ and so
$T$ are fixed there, but $F$ may run over those elements of
the orthonormal basis $\left\{u_{j,1/2}:\;j\ge 0\right\}$ having $
L$-eigenvalue
$1$. Then Shim$F$ runs over an orthogonal basis of even
Hecke normalized Maass-Hecke cusp forms of weight $0$
for ${\rm S}{\rm L}(2,{\rm Z})$, see [B-M], Theorem 1.2 and the last lines of
p. 982 of [D-I-T]. We see the coincidence with the
cuspidal sum of [Kw] in the above-mentioned special case. Note
that in the special case $h^{{\rm h}{\rm o}{\rm l}}\left(k\right)
=0$ the cuspidal sum of Theorem 3.1 of
[H-K] also has this form.

{\it Proof of Proposition 1.1.\/} We apply the Theorem of [B1]
for this $U$ and for $ $
$$u:={{\hbox{\rm Shim$F$}}\over {\sqrt {\left(\hbox{\rm Shim$F$},\hbox{\rm Shim$
F$}\right)_1}}}.\eqno (1.12)$$
Theorem 1.2 of [B-M] implies that we have a one-element sum in
the Theorem of [B1]. Then we get that
$$\left|\rho_u(1)\right|^2\left|\int_{D_1}\left|U(z)\right|^2u(z)
d\mu_z\right|^2=c_1\left|\rho_U\left(1\right)\right|^2\left|\rho_
F\left(1\right)\right|^2\left|\left(B_0\kappa\left(\overline U\right
),F\right)_4\right|^2,\eqno (1.13)$$
where $c_1>0$ is an absolute constant and $\rho_u(1)$, $\rho_U(1)$ are
the first Fourier coefficients of $u$ and $U$, respectively.

Formula (0.19) of [K-S] shows that
$$\left|\rho_F(1)\right|^2=c_2\left|\Gamma\left({{{1\over 2}+it}\over
2}\right)\right|^2\left|\rho_u(1)\right|^2L\left({1\over 2},\hbox{\rm Shim$
F$}\right),\eqno (1.14)$$
where $c_2>0$ is an absolute constant and $L\left(s,\hbox{\rm Shim$
F$}\right)$ is
the Hecke $L$-function of Shim$F$. We applied again
Theorem 1.2 of [B-M] to see that we have a one-element
sum in  [K-S], (0.19). We used also $t=2r$ and that (1.12) implies
$\left|\rho_u(1)\right|^2={1\over {\left(\hbox{\rm Shim$F$},\hbox{\rm Shim$
F$}\right)_1}}$.

By (2.4) of [B-K], which is a consequence of Watson's
identity (proved in [Wa], Theorem
3) we have that
$$\left|\int_{D_1}\left|U(z)\right|^2u(z)d\mu_z\right|^2=c_3{{\left
|\rho_U\left(1\right)\right|^4L\left({1\over 2},\hbox{\rm Shim$F$}\right
)L\left({1\over 2},\hbox{\rm Shim$F$}\otimes\hbox{\rm $\hbox{\rm sym}^
2U$}\right)}\over {\left(\hbox{\rm Shim$F$},\hbox{\rm Shim$F$}\right
)_1}}G_T\left(t\right)\eqno (1.15)$$
with
$$G_T\left(t\right):=\left|\Gamma\left({{{1\over 2}+it}\over 2}\right
)\right|^4\left|\Gamma\left({{{1\over 2}+2iT\pm it}\over 2}\right
)\right|^2,\eqno (1.16)$$
where $c_3>0$ is an absolute constant. We used in (2.4) of
[B-K] that the expressions
$$\left|\rho_U\left(1\right)\right|^2\left|\Gamma\left({1\over 2}
+iT\right)\right|^2L\left(1,\hbox{\rm $\hbox{\rm sym}^2U$}\right)
,\;\;\;\;{{L\left(1,\hbox{\rm $\hbox{\rm sym}^2u$}\right)\left|\Gamma\left
({1\over 2}+it\right)\right|^2}\over {\left(\hbox{\rm Shim$F$},\hbox{\rm Shim$
F$}\right)_1}}$$
are absolute constants, see [I-K], (5.101).

By (1.13), (1.14), (1.15), (1.16) and the fact that $L\left({1\over
2},\hbox{\rm Shim$F$}\right)\neq 0$
by our assumption $\rho_F(1)\neq 0$, and by (1.14) we get the statement.

{\bf 1.5. Statement of the holomorphic theorem.} First we
need some further definitions.

We introduce the Maass operators
$$K_k:=(z-\overline z){{\partial}\over {\partial z}}+k=iy{{\partial}\over {
\partial x}}+y{{\partial}\over {\partial y}}+k,\quad L_k:=(\overline
z-z){{\partial}\over {\partial\overline z}}-k=-iy{{\partial}\over {
\partial x}}+y{{\partial}\over {\partial y}}-k.$$
 \xrdef {MaassK} \xrdef {MaassL}
We will give the basic properties of these operators in Lemma 2.1 below. We just mention here that if $f$ is a Maass form of weight $k$, then $K_{k/2}f$ and $L_{k/2}
f$ are Maass forms of
weight $k+2$ and $k-2$, respectively.

If $k\ge 1$ is an integer, let $S_{2k+{1\over 2}}$ be the space of
holomorphic cusp forms of weight $2k+{1\over 2}$ \xrdef{hol} with the
multiplier system $\nu$ for the group $\Gamma_0(4)$. Let $f_{k,1}
,f_{k,2},...,f_{k,s_k}$ be an orthonormal basis of
$S_{2k+{1\over 2}}$, and write $g_{k,j}(z):=$$\left(\hbox{\rm Im$\,z$}\right
)^{{1\over 4}+k}f_{k,j}(z).$ We note that $g_{k,j}$ \xrdef{gkj} is
a Maass cusp form for $\Gamma_0(4)$ of weight $2k+{1\over 2}$, and
$\Delta_{}$$_{2k+{1\over 2}}g_{k,j}=\left(k+{1\over 4}\right)\left
(k-{3\over 4}\right)g_{k,j}$ (this follows easily from Lemma 2.1 below,
parts (v) and (iii)).

Suppose $u$ is a cusp form of weight $0$ for $SL(2,{\bf Z})$ with
$\Delta_0u$$=s(s-1)u$. For each $n\ge 0$, define
$$\left(\kappa_n(u)\right)(z):={{\left(K_{n-1}K_{n-2}\ldots K_1K_0
u\right)(4z)}\over {\left(s\right)_n\left(1-s\right)_n}},\eqno (1
.17)$$
 \xrdef{kappa(n)}
where $\left(a\right)_n:={{\Gamma\left(a+n\right)}\over {\Gamma\left
(a\right)}}$. \xrdef{poch} It is easy to check that $\kappa_n(u)$ is a cusp form of
weight $2n$ for the group $\Gamma_0$$\left(4\right)$.

{\bf THEOREM 1.2.} {\it For every integer} $n\ge 1$ {\it we have
that}
$$\sum_{j=1}^{s_n}\left(B_0\kappa_n\left(\overline {u_2}\right),g_{
n,j}\right)_4\overline {\left(B_0\kappa_n\left(\overline {u_1}\right
),g_{n,j}\right)_4}$$
 {\it equals the sum of}
$${6\over {\pi^{1/2}}}{{\delta_{u_1,u_2}\Gamma\left(n\pm it_1\right
)}\over {\Gamma\left(2n-{1\over 2}\right)\left(s_1\right)_n\left(
1-s_1\right)_n}}$$
 {\it and}
$$-{{24\pi\Gamma\left(n\pm it_1\right)\Gamma\left({1\over 2}\pm i
t_1-n\right)}\over {\Gamma\left(2n-{1\over 2}\right)\Gamma\left({
1\over 2}\pm it_2\right)}}{1\over {2\pi i}}\int_{\left({1\over 2}\right
)}{{\zeta\left(2S\right)L\left(S\right)\Gamma\left(S\right)\Gamma\left
(1-S\right)N\left(S,i\left({1\over 4}-n\right)\right)}\over {\left
(2\pi\right)^{2S}}}dS.\eqno (1.18)$$

 {\it The integral in
(1.18) is absolutely convergent.}

{\bf REMARK 1.3.} This result was informally announced in our
paper [B2], see pp 353-354. We decided to prove in this paper
only the nonholomorphic case, i.e. Theorem 1.1. Theorem 1.2 can be proved very similarly to the nonholomorphic case.

{\bf 1.6. Outline of the proof of Theorem 1.1.}

We have to give an expression for
$$\sum_{j=1}^{\infty}\chi\left(T_j\right)\left(B_0\kappa\left(\overline {
u_2}\right),u_{j,{1\over 2}}\right)_4\overline {\left(B_0\kappa\left
(\overline {u_1}\right),u_{j,{1\over 2}}\right)_4}+\hbox{\rm Eisenstein part}\eqno
(1.19)$$
with a weight function $\chi$. We can choose an automorphic
kernel $K(z,w)$ such that (1.19) equals
$$\int_{D_4}\left(\int_{D_4}\overline {B_0(z)}u_1\left(4z\right)K
(z,w)d\mu_z\right)B_0\left(w\right)\overline {u_2\left(4w\right)}
d\mu_w.\eqno (1.20)$$
By unfolding the inner integral here can be written as
$$\int_{\bbb H}\overline {B_0(z)}u_1\left(4z\right)k(z,w)d\mu_z\eqno
(1.21)$$
with a kernel function $k$. We now use geodesic polar
coordinates around $w$, so we have to compute the
integral on noneuclidean circles around $w$. We can
determine the Fourier expansion of $u_1$ on such circles
using an important theorem of Fay, which is recorded
in the present paper in Lemma 2.2. We get in this way
that (1.21) equals
$$\sum_{n=0}^{\infty}a_n\overline {B_n\left(w\right)}\left(K_{n-1}
K_{n-2}\ldots K_1K_0u_1\right)(4w),\eqno (1.22)$$
where
$$B_n:={1\over {n!}}K_{(n-1)+{1\over 4}}\ldots K_{{5\over 4}}K_{{1\over
4}}B_0,$$
and the coefficients $a_n$ are explicitly determined in
terms of the weight function $\chi$ and the
Laplace-eigenvalue of $u_1$. Inserting (1.22) in place of the
inner integral in (1.20) we get a weighted sum of integrals
$$\int_{D_4}\overline {B_n\left(w\right)}\left(K_{n-1}K_{n-2}\ldots
K_1K_0u_1\right)(4w)B_0\left(w\right)\overline {u_2\left(4w\right
)}d\mu_w.$$
This is the inner product involving $B_0\overline {B_n}$ what was
mentioned already in Remark 1.1. We show that $B_0\overline {B_n}$ is
a linear combination of Eisenstein series. Since the
Fourier coefficients of $K_{n-1}K_{n-2}\ldots K_1K_0u_1$ can be
given explicitly in terms of the Fourier coefficients of $u_1$, so by unfolding we get an
expression which contains the Rankin-Selberg $L$-function of $u_1$
and $\overline {u_2}$. Many problems occur concernig convergence and
the determination of the involved special functions, but these are the main steps of the proof of the
theorem.

To make the convergence problems easier we will first
impoose a stronger condition on the weight functions $\chi$
than the condition assumed in the theorem. This condition
will be the following:

We say that a function $\chi$ satisfies condition $D$ \xrdef {Cond(D)} if $\chi$ is
an even entire function satisfying that for every fixed $A,B>0$
the function $\left|\chi (z)\right|e^{\left|z\right|A}$ is bounded on the strip $\left
|\hbox{\rm Im}\,z\right|\le B$.

If a function $\chi$ satisfies Condition $D$, then it clearly
satisfies Condition $C_{\beta}$ for every $\beta >0$. Indeed, Condition $
D$ requires that $\chi$ decays faster than
exponentially on horizontal strips, while Condition $C_{\beta}$
allows exponential growth of a certain rate. We will first
prove the theorem for $\chi$ satisfying  Condition $D$. Then we
will show that it is relatively easy to extend the
statement for functions satisfying  $C_{\beta}$ with a suitable
$\beta >0$.

{\bf 1.7. Structure of the paper.} In Section 2 we list the
necessary notations and facts on automorphic functions.
In Section 3 we define the many types of special
functions occurring in the paper, give their properties
and prove some necessary lemmas on special functions.
We prove some very important lemmas needed for the
proof of Theorem 1.1 in Section 4, and we prove Theorem
1.1 in Section 5. However, the proofs of some important
lemmas on the kernel function and on the integral
transform are postponed to Section 6. We refer to the
statements of these lemmas in Section 5.

\noindent{\bf 2. Automorphic preliminaries}
\medskip

{\bf 2.1. Basic properties of the Maass operators.}

{\bf LEMMA 2.1.} {\it Let} $k,k_1,k_2\in {\bf R}$, $z\in\bbb H$, $
\gamma\in SL(2,{\bf Z})${\it , and let }
$f,g:\bbb H\rightarrow {\bf C}$ {\it be smooth functions. Then we have the following statements.}

{\it (i)} $K_{k_1-k_2}\left(f\overline g\right)=\left(K_{k_1}f\right
)\overline g+fK_{-k_2}\left(\overline g\right)$.

{\it (ii)} $\overline {\left(K_{-k}\overline f\right)}=L_kf$.

{\it (iii)} $\Delta_{2k}=L_{k+1}K_k+k$$\left(1+k\right)=K_{k-1}L_
k+k\left(k-1\right)$, $\Delta_{2k}L_{k+1}=L_{k+1}\Delta_{2k+2}$.

{\it (iv) If} $\Delta_{2k}f=s\left(s-1\right)f${\it , then for every} $
n\ge 0$ {\it we have }
$L_{k+1}\ldots L_{k+n}K_{k+n-1}\ldots K_kf=\left(-1\right)^n\left
(s+k\right)_n\left(1-s+k\right)_nf.$

{\it (v)} $f$ {\it is holomorphic if and only if} $K_k\left(y^{-k}\overline
f\right)=\overline {L_{-k}\left(y^{-k}f\right)}=0$.

{\it (vi)} $K_k\left(f\left(\gamma z\right)\left({{j_{\gamma}\left
(z\right)}\over {\left|j_{\gamma}\left(z\right)\right|}}\right)^{
-2k}\right)=\left({{j_{\gamma}\left(z\right)}\over {\left|j_{\gamma}\left
(z\right)\right|}}\right)^{-2k-2}\left(K_kf\right)\left(\gamma z\right
)$.

{\it Proof.\/} Parts (i), (ii) and (iii) follow by easy computations
using the definitions, and part (iv) follows easily from
(iii). Statement (iii) and (iv) are mentioned in [F],
formulas (6), (7) and (8). Part (v) is proved in Lemma 3.2 of
[R], and part (vi) is proved in Lemma 3.1 of [R]. The
proof is complete.

{\bf 2.2. Fourier expansions.} We first define the Fourier coefficients of Maass forms.
To do that the Whittaker functions $W_{\alpha ,\beta}$ are needed.
Their definition will be given in Section 3.7.

The three cusps for $\Gamma_0(4)$ are $\infty$, $0$ and $-{1\over
2}$. If $\eufm{a}$ denotes
one of these cusps, we take a scaling matrix $\sigma_{\eufm{a}}\in
SL(2,{\bf R})$ \xrdef{scaling}
as it is explained on p. 42 of [I]. We can easily
see that one can take
$$\sigma_{\infty}:=\left(\matrix{1&0\cr
0&1\cr}
\right),\qquad\sigma_0:=\left(\matrix{0&{{-1}\over 2}\cr
2&0\cr}
\right),\qquad\sigma_{-{1\over 2}}:=\left(\matrix{-1&{{-1}\over 2}\cr
2&0\cr}
\right).$$
The only cusp for $SL(2,{\bf Z})$ is $\infty$, and, of course, we take
the identity matrix $\sigma_{\infty}$ for scaling matrix also in this case.

If $\eufm{a}$ is a cusp for $\Gamma =SL(2,{\bf Z})$ or $\Gamma =\Gamma_
0(4)$, we define $\chi_{\eufm{a}}$ by
$$\nu\left(\sigma_{\eufm{a}}\left(\matrix{1&1\cr
0&1\cr}
\right)\sigma_{\eufm{a}}^{-1}\right)=e(-\chi_{\eufm{a}}),\qquad 0
\le\chi_{\eufm{a}}<1.$$
It is easy to check that $\chi_{\infty}=\chi_0=0$, and $\chi_{-{1\over
2}}={3\over 4}$. So
the cusps $0$ and $\infty$ are said to be singular, and $-1/2$ is
said to be nonsingular.

If $f$ is a Maass form of weight $l$, $\Delta_lf=s(s-1)f$ with some Re$
\, s\ge{1\over 2}$, $s={1\over 2}+it,$ and
$\eufm{a}$ is a cusp of $\Gamma$, then $f(\sigma_{\eufm{a}}z)\left
({{j_{\sigma_{\eufm{a}}}(z)}\over {\left|j_{\sigma_{\eufm{a}}}(z)\right
|}}\right)^{-l}$has the Fourier expansion \xrdef{rho(f,a,m)}
$$c_{f,\eufm{a}}(y)+\sum_{\matrix{m\in {\bf Z}\cr
m-\chi_{\eufm{a}}\neq 0\cr}
}\rho_{f,\eufm{a}}(m)W_{{l\over 2}{\rm s}{\rm g}{\rm n}\left(m-\chi_{\eufm{
a}}\right),it}\left(4\pi\left|m-\chi_{\eufm{a}}\right|y\right)e\left
(\left(m-\chi_{\eufm{a}}\right)x\right)$$
for $z=x+iy\in\bbb H$, and $c_{f,\eufm{a}}(y)=0$ if $\chi_{\eufm{
a}}\neq 0$, while it is a
linear combination of $y^s$ and $y^{1-s}$ for $s\neq{1\over 2}$ and of $
y^{1/2}$
and $y^{1/2}$$\log y$ for $s={1\over 2}$, if $\chi_{\eufm{a}}=0.$

We will need another type of Fourier expansion, namely Fourier expansion of Laplace-eigenfunctions on noneuclidean
circles. We reproduce here a theorem of Fay, which will
be important in the present paper. To state this theorem
we need geodesic polar coordinates: if $z_0\in\bbb H$ is fixed, then for every $
z\in\bbb H$ we
can uniquely write
$${{z-z_0}\over {z-\overline {z_0}}}=\tanh({r\over 2})e^{i\phi}\eqno
(2.1)$$
with $r>0$ and $0\le\phi <2\pi$. The invariant measure is expressed in these new coordinates as $
d\mu_z=\sinh rdrd\phi .$

{\bf LEMMA 2.2.} {\it Let} $k\in {\bf R}${\it ,} $s\in {\bf C}${\it , and let} $
f$ {\it be a smooth  function on} $H$
{\it satisfying} $\Delta_{2k}f=s\left(s-1\right )f${\it . If}
$z_0\in\bbb H$ {\it is given, then for every} $z\in\bbb H$ {\it we have
the absolutely } {\it convergent expansion}
$$f(z)\left({{z-\overline {z_0}}\over {z_0-\overline z}}\right)^k
=\sum_{n=-\infty}^{\infty}\left(f\right)_n(z_0)P_{s,k}^n(z,z_0)e^{
in\phi},$$
 {\it where} $r=r\left(z,z_0\right)>0$ {\it
and} $0\le\phi =\phi\left(z,z_0\right)<2\pi$ {\it are determined from} $
z$ {\it by (2.1), and}
$$P_{s,k}^n(z,z_0):=\left(\tanh({r\over 2})\right)^{\left|n\right|}\left
(1-\tanh^2({r\over 2})\right)^{k_n}F\left(s-k_n,1-s-k_n,1+\left|n\right
|,-y\right)$$
 {\it with} $y:={{\tanh^2({r\over 2})}\over {1-\tanh^2({r\over 2})}}${\it ,} $
k_n:=k{n\over {\left|n\right|}}$
{\it for} $n\neq 0${\it ,} $k_0:=\pm k${\it ,}
$$n!\left(f\right)_n(z_0):=\left(K_{k+n-1}\ldots K_{k+1}K_kf\right
)\left(z_0\right)\hbox{\rm \ for $n\ge 0$},$$
$$\left(-n\right)!\left(f\right)_n(z_0):=\overline {\left(K_{-k-n-
1}\ldots K_{1-k}K_{-k}\overline f\right)}\left(z_0\right)=\left(
L_{k+n+1}\ldots L_{k-1}L_kf\right)\left(z_0\right)\hbox{\rm \ for
$ n\le 0$}.$$
This follows from Theorems 1.1 and 1.2 of [F]. Lemma 2.2
was stated also in [B3], see Lemma 3.4 there. It is
explained there how to deduce Lemma 2.2 from the
theorems of Fay.

{\bf 2.3. The functions} $B_n${\bf .} If $z\in\bbb H$ is arbitrary, let $
T_z\in PSL(2,{\bf R})$ \xrdef{Tz} be such that $T_z$ is an upper triangular
matrix and $T_zi=z$. It is clear that $T_z$ is uniquely determined
by $z$, for $z=x+iy$ we have explicitly
$$T_z=\left(\matrix{y^{{1\over 2}}&xy^{{{-1}\over 2}}\cr
0&y^{{{-1}\over 2}}\cr}
\right).$$
If $z\in\bbb H$ is fixed, the function $\left(\hbox{\rm Im$\,z$}\right
)^{{1\over 4}}\theta\left(T_z\left(i{{1+L}\over {1-L}}\right)\right
)\left(1-L\right)^{-{1\over 2}}$ is holomorphic for
$\left|L\right|<1$, so it has a Taylor expansion
$$\hbox{\rm $\left(\hbox{\rm Im$\,z$}\right)^{{1\over 4}}\theta\left
(T_z\left(i{{1+L}\over {1-L}}\right)\right)\left(1-L\right)^{-{1\over
2}}$}=\sum_{n=0}^{\infty}B_n(z)L^n.\eqno (2.2)$$
 \xrdef {Bn}
We defined in this way a function $B_n(z)$ $(z\in\bbb H)$ for every $
n\ge 0$. For
$n=0$ this is in accordance with (1.2). These functions
satisfy also
$${1\over {n+1}}K_{n+{1\over 4}}B_n=B_{n+1}\eqno (2.3)$$
for every $n\ge 0$, this is proved in [B3], Lemma 6.1.
Indeed, this follows at once from (6.2) of [B3]. Formula
(2.3) implies that $B_n$ is a Maass form of weight $2n+{1\over 2}$ for $
\Gamma_0(4)$ and it has an additional
transformation formula
$$B_n\left({{-1}\over {4z}}\right)=e\left({{-1}\over 8}\right)\left
({z\over {\left|z\right|}}\right)^{{1\over 2}+2n}B_n(z)\eqno (2.4
)$$
for every $z\in\bbb H$, see (6.3), (6.4) and (6.5) of [B3].
These statements follow by induction using (1.3), (1.4), (2.3)
and Lemma 2.1 (vi).

{\bf 2.4. Rankin-Selberg} $L${\bf -functions.} It is known
that we have the functional equation
$${{\zeta\left(2S\right)L\left(S\right)}\over {\pi^{2S}}}\Gamma\left
({{S\pm it_1\pm it_2}\over 2}\right)={{\zeta\left(2\left(1-S\right
)\right)L\left(1-S\right)}\over {\pi^{2\left(1-S\right)}}}\Gamma\left
({{1-S\pm it_1\pm it_2}\over 2}\right)$$
for the Rankin-Selberg $L$-function defined in (1.6). We see
from this functional equation that the function $\zeta\left(2S\right
)L\left(S\right)$ is regular for
$S\neq 1$, and it has at most polynomial growth in vertical strips.

{\bf 2.5. Further notations.} We now explicitly give closures of fundamental domains of the
quotients $SL(2,{\bf Z})\setminus\bbb H$ and $\Gamma_0(4)\setminus\bbb
H$.

Let $D_1$ denote the closure of the standard fundamental
domain of the quotient $SL(2,{\bf Z})\setminus\bbb H$:
$$D_1:=\left\{z\in {\bf C}:\;\hbox{\rm Im$\, z>0$},\>-{1\over 2}\le\hbox{\rm Re$
\, z$}\le{1\over 2},\>\left|z\right|\ge 1\right\}.$$
 \xrdef{D1konk}
It is easy to check that the following set is a \xrdef{D4konk}
closure of a fundamental domain of  $\Gamma_0(4)\setminus\bbb H$:
$$D_4:=\bigcup_{j=0}^5\gamma_jD_1,$$
where
$$\gamma_j:=\left(\matrix{0&-1\cr
1&j\cr}
\right)\qquad (0\le j\le 3),\qquad\gamma_4:=\left(\matrix{1&0\cr
0&1\cr}
\right),\qquad\gamma_5:=\left(\matrix{1&0\cr
-2&1\cr}
\right).$$
 \xrdef {gammaj}
We always integrate over these fixed sets $D_1$ and $D_4$ in the sequel.

We denote by $R_l(D_4)$ \xrdef{Rl(D4)} the space of such smooth
automorphic functions of weight $l$$ $ for $\Gamma_0(4)$ for
which we have that for any integers $A,B,C\ge 0$ the function
$$\left(\max_{\eufm{a}}\hbox{\rm Im$\,\sigma_{\eufm{a}}^{-1}z$}\right
)^A\left|\left({{\partial^B}\over {\partial x^B}}{{\partial^C}\over {
\partial y^C}}f\right)(z)\right|$$
is bounded on $D_4$ (i.e. every partial derivative decays
faster than polynomially near each cusp on the fixed
fundamental domain $D_4$).

For $z,w\in\bbb H$, let
$$H(z,w):=i^{{1\over 2}}\left({{\left|z-\overline w\right|}\over {
(z-\overline w)}}\right)^{{1\over 2}},\eqno (2.5)$$
 \xrdef {H(z,w)}
as on p. 349 of [H]. It is easy to see that for any $T\in SL(2,{\bf R}
)$ we have
$${{H^2(Tz,Tw)}\over {H^2(z,w)}}=\left({{j_T(z)}\over {\left|j_T(
z)\right|}}\right)\left({{j_T(w)}\over {\left|j_T(w)\right|}}\right
)^{-1},$$
so
$${{H(Tz,Tw)}\over {H(z,w)}}=\left({{j_T(z)}\over {\left|j_T(z)\right
|}}\right)^{{1\over 2}}\left({{j_T(w)}\over {\left|j_T(w)\right|}}\right
)^{-{1\over 2}},\eqno (2.6)$$
since both sides lie in the right half-plane. Observe also that
$$H(w,z)=\overline {H(z,w)}.\eqno (2.7)$$

We now give the definition of the Eisenstein series of
weight $1/2$. For $\gamma_1,\gamma_2\in SL(2,{\bf R})$, we define
$$w(\gamma_1,\gamma_2):=j_{\gamma_1}(\gamma_2z)^{1/2}j_{\gamma_2}(
z)^{1/2}j_{\gamma_1\gamma_2}(z)^{-1/2},$$
the right-hand side is independent of $z\in\bbb H$. Clearly
$w=\pm 1$. For $\eufm{a}=0,\infty$, Re$\, s>1$, $z\in\bbb H$, define \xrdef {Eisenstein}
$$E_{\eufm{a}}\left(z,s,{1\over 2}\right):=\sum_{\gamma\in\Gamma_{\eufm{
a}}\setminus\Gamma_0(4)}\overline {\nu (\gamma )w\left(\sigma_{\eufm{
a}}^{-1},\gamma\right)}\hbox{\rm $($Im$\,\sigma_{\eufm{a}}^{-1}\gamma
z)^s$}\left({{j_{\sigma_{\eufm{a}}^{-1}\gamma}(z)}\over {\left|j_{
\sigma_{\eufm{a}}^{-1}\gamma}(z)\right|}}\right)^{-{1\over 2}},$$
where $\Gamma_{\eufm{a}}$ denotes the stability group of $\eufm{a}$
in $\Gamma_0(4)$.

Finally, we will use the notation $\Gamma_{\infty}:=\left\{\gamma\in SL(2
,{\bf Z}):\hbox{\rm \ $\gamma\infty =\infty$}\right\}$. \xrdef {GammaInfty}
The stability group of $\infty$ is clearly the same in $\Gamma_0(
4)$ and
$SL(2,{\bf Z})$.

\noindent{\bf 3. Preliminaries on special functions}
\medskip

{\bf 3.1. Generalized hypergeometric functions.} We define
these functions in the usual way: \xrdef {genhyp}
$$\hbox{\rm $_{q+1}$$F_q\left(\matrix{a_1,\ldots ,a_{q+1}\cr
b_1,\ldots ,b_q\cr}
;z\right):=$}\sum_{n=0}^{\infty}{{\left(a_1\right)_n\ldots\left(a_{
q+1}\right)_n}\over {n!\left(b_1\right)_n\ldots\left(b_q\right)_n}}
z^n.$$
Here the $b_i$ are not nonpsitive integers. We have absolute
convergence for $\left|z\right|<1$. The series is also absolutely
convergent for $\left|z\right|\le 1$ if we assume that
$\hbox{\rm Re}\, \left(\sum b_i-\sum a_i\right)>0$. $ $We will also use the notation $
F\left(\alpha ,\beta ,\gamma ;z\right)$ \xrdef{hyp} in place of
$_2$$F_1\left(\matrix{\alpha ,\beta\cr
\gamma\cr}
;z\right)$.

{\bf 3.2. Properties of} $_2$$F_1$ {\bf functions.} For $\hbox{\rm Re$
\, \alpha$}$, Re$\, \beta$, Re$\, \gamma >0$,
$-1<z<0$ and $-\hbox{\rm Re}\, \alpha$$,-\hbox{\rm Re}\, \beta <\sigma
<0$ we see by [S], (1.6.1.6) the Barnes-type integral
$$F\left(\alpha ,\beta ,\gamma ;z\right)={{\Gamma\left(\gamma\right
)}\over {\Gamma\left(\alpha\right)\Gamma\left(\beta\right)}}{1\over {
2\pi i}}\int_{(\sigma )}{{\Gamma\left(\alpha +s\right)\Gamma\left
(\beta +s\right)\Gamma\left(-s\right)}\over {\Gamma\left(\gamma +
s\right)}}\left(-z\right)^sds.\eqno (3.1)$$
This shows that $F\left(\alpha ,\beta ,\gamma ;z\right)$ extends analytically for $
z\notin [1,\infty )$.

For $\hbox{\rm Re$\, s$}<0$, $\hbox{\rm Re}\, \left(\alpha +s\right)>0$, $\hbox{\rm Re}\, \left
(\beta +s\right)>0$ we have that
$$\int_0^{\infty}x^{-s-1}F\left(\alpha ,\beta ,\gamma ;-x\right)d
x={{\Gamma\left(\gamma\right)\Gamma\left(\alpha +s\right)\Gamma\left
(\beta +s\right)\Gamma\left(-s\right)}\over {\Gamma\left(\alpha\right
)\Gamma\left(\beta\right)\Gamma\left(\gamma +s\right)}},\eqno (3.
2)$$ see [G-R],
p. 806, 7.511. For Re$\, \gamma >\hbox{\rm Re$\, \beta$}>0$, $z\notin [1,\infty )$ and any $
\alpha$ we have that
$$F\left(\alpha ,\beta ,\gamma ;z\right)={{\Gamma\left(\gamma\right
)}\over {\Gamma\left(\beta\right)\Gamma\left(\gamma -\beta\right)}}
\int_0^1t^{\beta -1}\left(1-t\right)^{\gamma -\beta -1}\left(1-tz\right
)^{-\alpha}dt,\eqno (3.3)$$ see [G-R], p. 995, 9.111.
If Re$\, \gamma >0$, $z\notin [0,\infty )$, and $\alpha$ and $\beta$ are any complex
numbers satisfying that $\alpha -\beta$ is not an integer, we have
$\hbox{\rm (idem}\left(\alpha ,\beta\right)$ means the same expression with $
\alpha$ and $\beta$
interchanged) that
$$F\left(\alpha ,\beta ,\gamma ;z\right)={{\Gamma\left(\gamma\right
)\Gamma\left(\beta -\alpha\right)}\over {\Gamma\left(\beta\right)
\Gamma\left(\gamma -\alpha\right)}}\left(-z\right)^{-\alpha}F\left
(\alpha ,\alpha +1-\gamma ,\alpha +1-\beta ;{1\over z}\right)+\hbox{\rm idem}\left
(\alpha ,\beta\right),\eqno (3.4)$$

see [S], (1.8.1.11). For $z\notin [1,\infty )$ we have by [G-R], p. 998, 9.131.1 that
$$F\left(\alpha ,\beta ,\gamma ;z\right)=\left(1-z\right)^{-\alpha}
F\left(\alpha ,\gamma -\beta ,\gamma ;{z\over {z-1}}\right)=\left
(1-z\right)^{\gamma -\alpha -\beta}F\left(\gamma -\alpha ,\gamma
-\beta ,\gamma ;z\right).\eqno (3.5)$$
{\bf 3.3. Properties of} $_3$$F_2$ {\bf functions.} Let $a$, $b$, $
c$ be such that $\hbox{\rm Re$\, a$}$, Re$\, b$, Re$\, c>0$ and the set $\left
\{a,b,c\right\}$ is symmetric with respect to
the real axis. We fix three numbers satisfying these
conditions throughout this subsection.

If $n$ is a nonnegative integer, the continuous dual Hahn polynomials are defined by \xrdef{Hahn}
$$S_n\left(x^2\right)=S_n\left(x^2;a,b,c\right):=\left(a+b\right)_
n\left(a+c\right)_n\,_3F_2\left(\matrix{-n,a+ix,a-ix\cr
a+b,a+c\cr}
;1\right),\eqno (3.6)$$
see [A-A-R], (6.10.2). Formula (3.6) is symmetric in
the parameters $a$, $b$, $c$, this follows from the identity
$$_3F_2\left(\matrix{A,B,C\cr
D,E\cr}
;1\right)={{\Gamma\left(E\right)\Gamma\left(D+E-A-B-C\right)}\over {
\Gamma\left(E-A\right)\Gamma\left(D+E-B-C\right)}}\,_3F_2\left(\matrix{
A,D-B,D-C\cr
D,D+E-B-C\cr}
;1\right)$$
valid for $\hbox{\rm Re}\, \left(D+E-B-C-A\right)>0,$ $\hbox{\rm Re}\, \left
(E-A\right)>0,$ see
Corollary 3.3.5 of [A-A-R]. These polynomials form a
complete orthogonal system in $L^2$$\left(\left(0,\infty\right),w_{
a,b,c}\left(x\right)dx\right)$
with the weight function \xrdef {weight}
$$w_{a,b,c}\left(x\right):={1\over {2\pi}}{{\Gamma\left(a\pm ix\right
)\Gamma\left(b\pm ix\right)\Gamma\left(c\pm ix\right)}\over {\Gamma\left
(\pm 2ix\right)}}.\eqno (3.7)$$
Indeed, we have the relations
$$\int_0^{\infty}w_{a,b,c}\left(x\right)S_m\left(x^2\right)S_n\left
(x^2\right)dx=\delta_{mn}\Gamma\left(n+a+b\right)\Gamma\left(n+a+
c\right)\Gamma\left(n+b+c\right)n!,\eqno (3.8)$$
where $\delta_{mn}$ is the Kronecker delta symbol, see [A-A-R],
(6.10.7). Completeness of the system follows from
Theorem 6.5.2 of [A-A-R], taking into account that
$w_{a,b,c}\left(x\right)$ decays exponentially as $x\rightarrow +
\infty$.

We can deduce a pointwise upper bound from (3.8). This
bound is weak, but it will be enough for our purposes.

{\bf LEMMA 3.1.} {\it There is a positive} $M$ {\it such that}
$$\left|{{S_n\left(x^2;a,b,c\right)}\over {\left(a+b\right)_n\left
(a+c\right)_n}}\right|\le Me^{{{\pi}\over 2}\left|x\right|}\left(
1+n\right)^M$$
 {\it for every
integer} $n\ge 0$ {\it and every real} $x$.

{\it Proof.\/} It is enough to show that for every real $y$ we
have
$$\max_{y\le x\le y+1}\left|{{S_n\left(x^2;a,b,c\right)}\over {\left
(a+b\right)_n\left(a+c\right)_n}}\right|\le Me^{{{\pi}\over 2}\left
|y\right|}\left(1+n\right)^M\eqno (3.9)$$
with a suitable $M$. The classical Markov inequality states that for any
polynomial $p$ of degree $n$ we have
$$\max_{-1\le t\le 1}\left|p^{(1)}\left(t\right)\right|\le n^2\max_{
-1\le t\le 1}\left|p\left(t\right)\right|.$$
This is proved e.g. in [B-E], Theorem 5.1.8. Then we see
that if the left-hand side of (3.9) is $m$, then there is a subinterval $
I$ of
$[y,y+1]$ such that the length of $I$ is $\gg n^{-2}$, and
$\left|{{S_n\left(x^2;a,b,c\right)}\over {\left(a+b\right)_n\left
(a+c\right)_n}}\right|\gg m$ for every $x\in I$. Then we get the
lemma by (3.7), the $m=n$ case of (3.8) and the Stirling
formula.

{\bf LEMMA 3.2.} {\it If}  $\hbox{\rm $\hbox{\rm Re}\, \gamma$}$, $\hbox{\rm Re}
\, A>0$ {\it and} $\hbox{\rm Re}\, B$ {\it is large enough in terms of} $
a,b,c$
{\it and} $\hbox{\rm Re$\, \gamma$}${\it , then we have for every real} $
x$ {\it that}
$$\sum_{n=0}^{\infty}{{_3F_2\left(\matrix{-n,a+ix,a-ix\cr
a+b,a+c\cr}
;1\right)}\over {n!}}{{\left(\gamma\right)_n\left(A\right)_n}\over {\left
(A+B\right)_n}}$$
{\it equals}
$${{\Gamma\left(a+b\right)\Gamma\left(a+c\right)\Gamma\left(A+B\right
)}\over {\Gamma\left(\gamma\right)\Gamma\left(A\right)\Gamma\left
(B\right)\Gamma\left(a\pm ix\right)\Gamma\left(A+B-\gamma\right)}}$$
{\it times}
$${1\over {2\pi i}}\int_{\left(-C\right)}{{\Gamma\left(\gamma +s\right
)\Gamma\left(A+s\right)\Gamma\left(a\pm ix+s\right)\Gamma\left(-s\right
)\Gamma\left(B-\gamma -s\right)}\over {\Gamma\left(a+b+s\right)\Gamma\left
(a+c+s\right)}}ds$$
{\it with} $0<C<\min\left(\hbox{\rm Re}\, a,\hbox{\rm Re}\, \gamma ,\hbox{\rm Re}\,
A\right)$.

{\it Proof.\/} For any complex $\gamma$ and real $x$ for $\left|t\right
|<1/2$
we have
$$\sum_{n=0}^{\infty}{{_3F_2\left(\matrix{-n,a+ix,a-ix\cr
a+b,a+c\cr}
;1\right)}\over {n!}}\left(\gamma\right)_nt^n=\left(1-t\right)^{-
\gamma}\,_3F_2\left(\matrix{\gamma ,a+ix,a-ix\cr
a+b,a+c\cr}
;{t\over {t-1}}\right).\eqno (3.10)$$
Indeed, this follows easily by inserting on the left-hand side the defining series of $_
3F_2$, changing the
summations, and using for every nonnegative integer $k$ that
$$\sum_{n=k}^{\infty}{{\left(-n\right)_k\left(\gamma\right)_nt^n}\over {
n!}}=\left(\gamma\right)_k\left(1-t\right)^{-\gamma}\left({t\over {
t-1}}\right)^k,$$
which follows from (3.5). For $0<t<{1\over 2}$ the right-hand side here equals
$${{\Gamma\left(a+b\right)\Gamma\left(a+c\right)}\over {\Gamma\left
(\gamma\right)\Gamma\left(a\pm ix\right)}}{{\left(1-t\right)^{-\gamma}}\over {
2\pi i}}\int_{\left(-C\right)}{{\Gamma\left(\gamma +s\right)\Gamma\left
(a\pm ix+s\right)\Gamma\left(-s\right)}\over {\Gamma\left(a+b+s\right
)\Gamma\left(a+c+s\right)}}\left({t\over {1-t}}\right)^sds,\eqno
(3.11)$$
which can be seen by shifting the integration to the
right. Then using Lemma 3.1 by analytic continuation we see that the
left-hand side of (3.10) equals (3.11)
for any $0<t<1$. Multiplying by $t^{A-1}$$\left(1-t\right)^{B-1}$, integrating from 0 to 1 and using (3.3) with $
z=0$ we obtain the lemma.

We have the difference equation
$$nS_n\left(x^2\right)=B\left(x\right)S_n\left(\left(x+i\right)^2\right
)-\left(B\left(x\right)+D\left(x\right)\right)S_n\left(x^2\right)
+D\left(x\right)S_n\left(\left(x-i\right)^2\right)\eqno (3.12)$$
for every $n\ge 0$, where we write
$$B\left(x\right)={{\left(a-ix\right)\left(b-ix\right)\left(c-ix\right
)}\over {\left(-2ix\right)\left(1-2ix\right)}},\;D\left(x\right)={{\left
(a+ix\right)\left(b+ix\right)\left(c+ix\right)}\over {\left(2ix\right
)\left(1+2ix\right)}},$$
see [A-A-R], (6.10.9). This relation has the following
consequence.

{\bf LEMMA 3.3.} {\it Let} $\chi$ {\it be a function satisfying Condition} $
D${\it . For any} $A>0$ {\it we have for integers} $n\ge 0$ {\it that }
$$C_{n,\chi}\left(a,b,c\right):=\int_{-\infty}^{\infty}\chi (x){{
S_n\left(x^2;a,b,c\right)}\over {\left(a+b\right)_n\left(a+c\right
)_n}}w_{a,b,c}\left(x\right)dx\ll_{\chi ,a,b,c,A}\left(1+n\right)^{
-A}.\eqno (3.13)$$
{\it Proof.\/} We substitute (3.12) into (3.13), and we shift the integration to
${\rm I}{\rm m}\,x=1$ in the case of $S_n\left((x-i)^2\right)$, and to
${\rm I}{\rm m}\,x=-1$ in the case of $S_n\left((x+i)^2\right)$. We do not cross any
pole, and we get for $nC_{n,\chi}\left(a,b,c\right)$ an expression of type
(3.13), but with a new function in place of $\chi$ satisfying
Condition $D$. These facts can be checked using (3.7). We iterate this step many times, and
then we apply Cauchy-Schwarz inequality and use (3.8)
with $m=n$. By the properties of $\chi$ this proves the
lemma.

{\bf 3.4. Some integral formulas.} For $0<\hbox{\rm Re$\, a_1$}$, Re$\,
a_2$, Re$\, b_1$, Re$\, b_2$, Re$\, b_3$, assuming that $b_4+a_1$
is not a nonpositive integer, we have that
$${1\over {2\pi i}}\int_{(0)}{{\Gamma\left(a_1-s\right)\Gamma\left
(a_2-s\right)\Gamma\left(b_1+s\right)\Gamma\left(b_2+s\right)\Gamma\left
(b_3+s\right)}\over {\Gamma\left(b_4+s\right)}}ds\eqno (3.14)$$
equals
$$_3F_2\left(\matrix{a_1+b_1,a_1+b_2,b_4-b_3\cr
a_1+a_2+b_1+b_2,b_4+a_1\cr}
;1\right){{\Gamma\left(a_1+b_3\right)\prod_{k=1}^2\left(\Gamma\left
(a_1+b_k\right)\Gamma\left(a_2+b_k\right)\right)}\over {\Gamma\left
(a_1+a_2+b_1+b_2\right)\Gamma\left(b_4+a_1\right)}},\eqno (3.15)$$
see [S], (4.2.2.1).

In the special case $b_4=a_1+a_2+b_1+b_2+b_3$ we have the
following statement. For $0<\hbox{\rm Re$\, a_1$}$, Re$\, a_2$, Re$\, b_1$, Re$\,
b_2$, Re$\, b_3$ we have that
$${1\over {2\pi i}}\int_{(0)}{{\Gamma\left(a_1-s\right)\Gamma\left
(a_2-s\right)\Gamma\left(b_1+s\right)\Gamma\left(b_2+s\right)\Gamma\left
(b_3+s\right)}\over {\Gamma\left(b_1+b_2+b_3+a_1+a_2+s\right)}}ds\eqno
(3.16)$$
equals
$${{\prod_{k=1}^3\left(\Gamma\left(a_1+b_k\right)\Gamma\left(a_2+
b_k\right)\right)}\over {\prod_{1\le k<l\le 3}\Gamma\left(b_k+b_l
+a_1+a_2\right)}}.\eqno (3.17)$$
This is the Second Barnes Lemma, see [S], (4.2.2.2).

{\bf LEMMA 3.4.} {\it For} $\hbox{\rm $0<\hbox{\rm Re$\, \alpha\le\hbox{\rm Re$\,
\beta$}$}$}$, $0<\hbox{\rm Re$\, a\le\hbox{\rm Re$\, b$}$}$, $0<\hbox{\rm Re$\,
\gamma$}$,
$0<\hbox{\rm Re$\, c<\hbox{\rm Re$\, a+$}$}\hbox{\rm Re$\, \alpha ,$ }${\it if} $
\gamma +a-c$ {\it is not a nonpositive integer, we have that }
$$\int_0^{\infty}\hbox{\rm $F\left(\alpha ,\beta ,\gamma ;-u\right
)u^{c-1}F\left(a,b,c;-u\right)du$}\eqno (3.18)$$
{\it equals the product of}
$${{\Gamma\left(\gamma\right)\Gamma\left(c\right)\Gamma\left(a-c+
\alpha\right)\Gamma\left(a-c+\beta\right)\Gamma\left(b-c+\alpha\right
)\Gamma\left(b-c+\beta\right)}\over {\Gamma\left(\alpha\right)\Gamma\left
(\beta\right)\Gamma\left(b\right)\Gamma\left(a+b-2c+\alpha +\beta\right
)\Gamma\left(\gamma +a-c\right)}}$$
{\it and}
$$_3F_2\left(\matrix{a-c+\alpha ,a-c+\beta ,\gamma -c\cr
a+b-2c+\alpha +\beta ,\gamma +a-c\cr}
;1\right).$$
{\it Proof.\/} Applying (3.1) for the first factor in (3.18) with $
\sigma$ satisfying
$$-\hbox{\rm Re$\, \alpha$},-\hbox{\rm Re$\, c$}<\sigma <\hbox{\rm Re$\, \left
(a-c\right)$},0$$
and then using (3.2) we get that (3.18) equals
$${{\Gamma\left(\gamma\right)\Gamma\left(c\right)}\over {\Gamma\left
(\alpha\right)\Gamma\left(\beta\right)\Gamma\left(a\right)\Gamma\left
(b\right)}}{1\over {2\pi i}}\int_{(\sigma )}{{\Gamma\left(a-c-s\right
)\Gamma\left(b-c-s\right)\Gamma\left(\alpha +s\right)\Gamma\left(
\beta +s\right)\Gamma\left(c+s\right)}\over {\Gamma\left(\gamma +
s\right)}}ds.$$
By the equality of (3.14) and (3.15) this gives the statement
of the lemma.

{\bf 3.5. A hypergeometric integral transform.} Our aim here is to prove Lemma 3.7 below, to prepare
its proof we need the identities proved in the next two lemmas.

{\bf LEMMA 3.5.} {\it For every integer} $n\ge 0$ {\it and for every real} $
t$ {\it we have that}
$$\int_0^{\infty}\hbox{\rm $F\left({3\over 4}+it,{3\over 4}-it,1;
-u\right){{1-\left({u\over {1+u}}\right)^{n+1}}\over {n+1}}du$}={{
\Gamma\left({3\over 4}\pm it\right)\,_3F_2\left(\matrix{{3\over 4}
+it,{3\over 4}-it,-n\cr
{3\over 2},2\cr}
;1\right)}\over {\Gamma\left({3\over 2}\right)}}.$$
{\it Proof.\/} Let $m\ge 0$ be an integer. Writing $b=c=m+1$, $\gamma
=1$
we get from Lemma 3.4 and from (3.5) that
$$\int_0^{\infty}\hbox{\rm $F\left(\alpha ,\beta ,1;-u\right){{u^
m}\over {\left(1+u\right)^a}}du$}=\sum_{k=0}^m{{\left(-m\right)_k
\Gamma\left(a-m-1+\alpha +k\right)\Gamma\left(a-m-1+\beta +k\right
)}\over {k!\Gamma\left(a-m-1+\alpha +\beta +k\right)\Gamma\left(a
-m+k\right)}}\eqno (3.19)$$
under the conditions $\hbox{\rm $0<\hbox{\rm Re$\, \alpha\le\hbox{\rm Re$\,
\beta$}$}$}$, $0<\hbox{\rm Re$\, a\le m+1$}$,
$\hbox{\rm $m+1-\hbox{\rm Re$\, a$}<\hbox{\rm Re$\, \alpha$}$}$, assuming that $
a$ is not an integer. We estimate the hypergeometric
function on the left-hand side by (3.1), and we see by
analytic continuation in $a$ that it
is enough to assume $0<\hbox{\rm Re$\, \alpha\le\hbox{\rm Re$\, \beta$}$}$ and $
m+1$$-\hbox{\rm Re$\, a$}<\hbox{\rm Re$\, \alpha$}$. So
assuming $1<\hbox{\rm Re$\, \alpha\le\hbox{\rm Re$\, \beta$}$}$ we see that (3.19) is true for
$a=m$. Taking the difference of (3.19) for $m=0$ and
$m=n+1$ we get the lemma by analytic continuation in $\alpha$ and $
\beta$.

{\bf LEMMA 3.6.} {\it For any integer} $n\ge 0$ {\it and for any} $
u>0$ {\it we have that}
$$\int_0^{\infty}F\left({3\over 4}-it,{3\over 4}+it,1,-u\right){{
\Gamma^2\left({3\over 4}\pm it\right)\Gamma\left({1\over 4}\pm it\right
)}\over {\Gamma\left(\pm 2it\right)}}\,_3F_2\left(\matrix{{3\over
4}+it,{3\over 4}-it,-n\cr
{3\over 2},2\cr}
;1\right)\,dt\eqno (3.20)$$
{\it equals}
$$2\pi\Gamma\left({3\over 2}\right)\left(1+u\right)^{-{1\over 2}}{{
1-\left({u\over {1+u}}\right)^{n+1}}\over {n+1}}.$$
{\it Proof.\/} First note that for any $\hbox{\rm Re$\, s={1\over 2}$ }$and any integer $
k\ge 0$
we have by (3.7) and (3.8) that
$$\int_0^{\infty}{{\Gamma\left({1\over 4}\pm it\right)\Gamma\left
({3\over 4}\pm it+k\right)\Gamma\left({3\over 4}\pm it+s\right)}\over {
\Gamma\left(\pm 2it\right)}}\,dt=2\pi\Gamma\left(1+k\right)\Gamma\left
(1+s\right)\Gamma\left({3\over 2}+k+s\right).$$
Then by (3.1) we have that (3.20) equals
$$2\pi\sum_{k=0}^n{{\left(-n\right)_k}\over {\left({3\over 2}\right
)_k\left(2\right)_k}}{1\over {2\pi i}}\int_{\left(-1/2\right)}\Gamma\left
(-s\right)\Gamma\left({3\over 2}+k+s\right)u^sds.$$
We can compute the integral here by (3.1) and (3.5), and
we get that (3.20) equals
$$2\pi\Gamma\left({3\over 2}\right)\sum_{k=0}^n{{\left(-n\right)_
k}\over {\left(2\right)_k}}\left(1+u\right)^{-{3\over 2}-k}.$$
The lemma follows by the binomial theorem.

The integral transform (3.21) below is a special case of
the so-called {\it Jacobi transform}, see e.g. [K]. The
inversion formula of this transform is also proved there
in Theorem 4.2, but since it is not hard to prove it using the results on
continuous dual Hahn polynomials mentioned above, so we include a proof.

{\bf LEMMA 3.7.} {\it Let} $\chi$ {\it be a function satisfying Condition} $
D${\it . For} $u\ge 0$ {\it define  }
$$k_{\chi}(u):={1\over {\pi}}(u+1)^{{1\over 4}}\int_0^{\infty}F\left
({3\over 4}-it,{3\over 4}+it,1,-u\right)\left|{{\Gamma\left({1\over
4}+it\right)\Gamma\left({3\over 4}+it\right)}\over {\Gamma\left(2
it\right)}}\right|^2\chi (t)dt.\eqno (3.21)$$
{\it Then the following statements hold.}

{\it (i) The function} $k_{\chi}^{(j)}(u)$$(u+1)^A$ {\it is bounded on} $
[0,\infty )$ {\it for every} $A>0$ {\it and} $j\ge 0$.

{\it (ii) For every real} $t$ {\it we have }
$$\chi (t)={1\over 2}\int_0^{\infty}(u+1)^{{1\over 4}}F\left({3\over
4}-it,{3\over 4}+it,1,-u\right)k_{\chi}(u)du.$$
{\it Proof.\/} By (3.4) we know for real $t$ that
$$F\left({3\over 4}-it,{3\over 4}+it,1,-u\right)\left|{{\Gamma\left
({1\over 4}+it\right)\Gamma\left({3\over 4}+it\right)}\over {\Gamma\left
(2it\right)}}\right|^2=\phi (u,t)+\phi (u,-t),$$
where
$$\phi (u,t)={{\Gamma\left({1\over 4}-it\right)\Gamma\left({3\over
4}-it\right)}\over {\Gamma\left(-2it\right)}}u^{it-{3\over 4}}F\left
({3\over 4}-it,{3\over 4}-it,1-2it,-{1\over u}\right),$$
hence
$$k_{\chi}(u)={1\over {\pi}}(u+1)^{{1\over 4}}\int_{-\infty}^{\infty}
\phi (u,t)\chi (t)dt.$$
Now, if $u$ is large, we push the line of integration
upwards to a line $\hbox{\rm Im$\,t=B$ }$with a large positive number
$B$ depending on $A$ and $j$, and using (3.3) we get (i). Indeed,
for small $u$ statement (i) is trivial, using the very definition of $
k_{\chi}(u)$ and (3.3).

Let $\left\{a,b,c\right\}=\left\{{3\over 4},{3\over 4},{5\over 4}\right
\}$, and let us write
$${{\chi (t)}\over {\Gamma\left({3\over 4}\pm it\right)}}=\sum_{n
=0}^{\infty}a_n\,{{S_n\left(t^2;{3\over 4},{3\over 4},{5\over 4}\right
)}\over {\left({3\over 2}\right)_n\left(2\right)_n}}\eqno (3.22)$$
in the space $L^2$$\left(\left(0,\infty\right),w_{a,b,c}\left(t\right
)dt\right)$. It follows from (3.8)
and Lemma 3.3 that $a_n\ll_{\chi ,A}\left(1+n\right)^{-A}$ for every $
A>0$.
Then Lemma 3.1 shows that the right-hand side of (3.22)
is a continuous function, so (3.22) is valid pointwise for
every $t>0$. We also see applying Lemma 3.1 and (3.1) that if we
express $\chi (t)$ from (3.22) and substitute the obtained
expressioon into (3.21), then we can integrate there term
by term. Then from (3.6) and Lemma 3.6 we get that
$$k_{\chi}(u)=2\Gamma\left({3\over 2}\right)\left(1+u\right)^{-{1\over
4}}\sum_{n=0}^{\infty}a_n\,{{1-\left({u\over {1+u}}\right)^{n+1}}\over {
n+1}}.$$
By Lemma 3.5, (3.6) and (3.22) we obtain (ii). The lemma is proved.

{\bf 3.6. Properties of} $_7$$F_6$ {\bf and} $_4$$F_3$ {\bf functions.} For complex
$A,B,C,D,E$ and $F$ satisfying $\hbox{\rm Re}\, \left(2+2A-B-C-D-E-F\right
)>0$
let \xrdef{w}
$$W\left(A;B,C,D,E,F\right):=\>_7F_6\left(\matrix{A,1+{A\over 2},
B,C,D,E,F\cr
{A\over 2},B^{\ast},C^{\ast},D^{\ast},E^{\ast},F^{\ast}\cr}
;1\right)$$
where $B^{\ast},C^{\ast},D^{\ast},E^{\ast},F^{\ast}$ are given by
$$B+B^{\ast}=C+C^{\ast}=D+D^{\ast}=E+E^{\ast}=F+F^{\ast}=1+A.$$
If $a,b,c,d$ are complex numbers, then writing
$$\tilde {a}:={1\over 2}\left(a+b+c+d-1\right),\qquad\tilde {b}:={1\over
2}\left(a+b-c-d+1\right),\eqno (3.23)$$
$$\tilde {c}:={1\over 2}\left(a-b+c-d+1\right),\qquad\tilde {d}:={1\over
2}\left(a-b-c+d+1\right),\eqno (3.24)$$
we define the Wilson function \xrdef{wilson}
$\phi_{\lambda}\left(x\right)=\phi_{\lambda}\left(x;a,b,c,d\right
)$ by the formula
$$\phi_{\lambda}\left(x\right):={{\Gamma\left(\tilde a+\tilde b+\tilde
c+i\lambda\right)W\left(\tilde a+\tilde b+\tilde c+i\lambda -1;a+
ix,a-ix,\tilde a+i\lambda ,\tilde b+i\lambda ,\tilde c+i\lambda\right
)}\over {\Gamma\left(a+b\right)\Gamma\left(a+c\right)\Gamma\left(
1+a-d\right)\Gamma\left(1-\tilde d-i\lambda\right)\Gamma\left(\tilde
b+c+i\lambda\pm ix\right)}},$$
as it was introduced by Groenevelt in [G1] in formula
(3.2). This definition is meaningful if the $\Gamma$-function is regular
at the point $\tilde {a}+\tilde {b}+\tilde {c}+i\lambda$ and
$\hbox{\rm Re}$$\, \left(1-\tilde d-i\lambda\right)>0$. However, $\phi_{
\lambda}\left(x\right)$ is an entire function in $ $$\left(x,\lambda\right
)\in {\bf C}^2$ (see [G1], below formula
(3.3)).

The Wilson function $\phi_{\lambda}\left(x;a,b,c,d\right)$ is symmetric in the
parameters $a$, $b$, $c$, $1-d$, see [G2], Lemma 5.3 (ii). We
have the symmetry relation
$$\phi_{\lambda}\left(x;a+i\omega ,b+iy,b-iy,1-a+i\omega\right)=\phi_{
\omega}\left(y;a+i\lambda ,b+ix,b-ix,1-a+i\lambda\right),\eqno (3.25)$$
see [G2], Lemma 5.3 (i).

We have the identity that
$${1\over {2\pi i}}\int_{-i\infty}^{i\infty}{{\Gamma\left(a\pm ix
+R\right)\Gamma\left(\tilde a\pm i\lambda +R\right)\Gamma\left(-R\right
)\Gamma\left(1-a-d-R\right)}\over {\Gamma\left(a+b+R\right)\Gamma\left
(a+c+R\right)}}dR\eqno (3.26)$$
equals
$$\Gamma\left(a\pm ix\right)\Gamma\left(\tilde a\pm i\lambda\right
)\Gamma\left(1-d\pm ix\right)\Gamma\left(1-\tilde d\pm i\lambda\right
)\phi_{\lambda}\left(x;a,b,c,d\right),\eqno (3.27)$$
assuming that $\Gamma$ is regular at the points $a\pm ix$, $\tilde {
a}\pm i\lambda$,
$1-d\pm ix$, $1-\tilde {d}\pm i\lambda$. Here the poles of the functions $
\Gamma\left(a\pm ix+R\right)$, $\Gamma\left(\tilde a\pm i\lambda
+R\right)$ lie to
the left of the path of integration, and the poles of the
functions $\Gamma\left(-R\right)$, $\Gamma\left(1-a-d-R\right)$
lie to the right of it, and $\tilde {a}$$,\tilde {b},\tilde {c},\tilde {
d}$ are defined above.
This can be seen by shifting the integration to the
right in (3.26) above and applying (3.3) of [G1].

We need some important identities, and in order to state
them we need further notations. For complex $A,B,C,D,E$
and $F$ \xrdef {psi} let
$$\psi\left(A;B,C,D,E,F\right):={{\Gamma\left(1+A\right)W\left(A;B
,C,D,E,F\right)}\over {\Gamma\left(B^{\ast}\right)\Gamma\left(C^{
\ast}\right)\Gamma\left(D^{\ast}\right)\Gamma\left(E^{\ast}\right
)\Gamma\left(F^{\ast}\right)\Gamma\left(1-S\right)}},\eqno (3.28)$$
where
$$S:=B+C+D+E+F-2A-1,\eqno (3.29)$$
see (2.1) of [W] or p 127 of [S].

We can check that we have
$$\phi_{\lambda}\left(x;a,b,c,d\right)=\psi\left(\tilde a+\tilde
b+\tilde c+i\lambda -1;a+ix,a-ix,\tilde a+i\lambda ,\tilde b+i\lambda
,\tilde c+i\lambda\right)\eqno (3.30)$$
with this notation.

Let us write
$$\,\,_4F_3^{\ast}\left(\matrix{A,B,C,D\cr
E,F,G\cr}
;1\right):={{\Gamma\left(A\right)\Gamma\left(B\right)\Gamma\left(
C\right)\Gamma\left(D\right)}\over {\Gamma\left(E\right)\Gamma\left
(F\right)\Gamma\left(G\right)}}\,_4F_3\left(\matrix{A,B,C,D\cr
E,F,G\cr}
;1\right).\eqno (3.31)$$
Then assuming $E+F+G-A-B-C-D=1$ and $\hbox{\rm Re}\, \left(1+A-G\right
)>0$ we have that
$$\,_4F_3^{\ast}\left(\matrix{A,B,C,D\cr
E,F,G\cr}
;1\right)-\,_4F_3^{\ast}\left(\matrix{1+A-G,1+B-G,1+C-G,1+D-G\cr
1+E-G,1+F-G,2-G\cr}
;1\right)=P_1P_2\eqno (3.32)$$
with the abbreviations
$$P_1:={{\Gamma\left(A\right)\Gamma\left(B\right)\Gamma\left(C\right
)\Gamma\left(D\right)\Gamma\left(1+A-G\right)\Gamma\left(1+B-G\right
)\Gamma\left(1+C-G\right)\Gamma\left(1+D-G\right)}\over {\Gamma\left
(G\right)\Gamma\left(1-G\right)}}\eqno (3.33)$$
and
$$P_2:=\psi\left(B+C+D-G;B,C,D,E-A,F-A\right).\eqno (3.34)$$
This follows by some computations from (2.4.4.3) of [S].
See also (2.3) of [W].

Assuming $\hbox{\rm Re}\, \left(2+2A-B-C-D-E-F\right)>0$ and $\hbox{\rm Re}\, \left
(B+D-A\right)>0$ we have that
$${{\psi\left(A;B,C,D,E,F\right)\sin\left(\pi\left(D+E+F-A\right)\right
)}\over {\Gamma\left(B+D-A\right)\Gamma\left(B+E-A\right)\Gamma\left
(B+F-A\right)\Gamma\left(1-C\right)}}\eqno (3.35)$$
equals the sum of
$${{\psi\left(E+F-C;E,F,1+A-B-C,1+A-C-D,E+F-A\right)\sin\left(\pi\left
(B-A\right)\right)}\over {\Gamma\left(1+A-E-F\right)\Gamma\left(1
+A-D-F\right)\Gamma\left(1+A-D-E\right)\Gamma\left(1-S\right)}}\eqno
(3.36)$$
and
$${{\psi\left(2B-A;B,B+C-A,B+D-A,B+E-A,B+F-A\right)\sin\left(\pi\left
(C-S\right)\right)}\over {\Gamma\left(D\right)\Gamma\left(E\right
)\Gamma\left(F\right)\Gamma\left(1+A-B-C\right)}},\eqno (3.37)$$
using the notation (3.29). This is formula (2.7) of [W] (see
also (4.3.7.8) of [S]).

{\bf 3.7. Whittaker functions.} For complex numbers $\alpha$, $\beta$ satisfying $
0<{1\over 2}-\hbox{\rm Re}\, \alpha -\left|\hbox{\rm Re}\, \beta\right|$ and for $
y>0$ we define the Whittaker function $W_{\alpha ,\beta}\left(y\right
)$ \xrdef{Whittaker} by
the formula
$$W_{\alpha ,\beta}\left(y\right):={{e^{-{y\over 2}}}\over {2\pi i}}
\int_{\left(\sigma\right)}{{\Gamma\left(v\right)\Gamma\left({1\over
2}-\alpha\pm\beta -v\right)}\over {\Gamma\left({1\over 2}-\alpha\pm
\beta\right)}}y^{\alpha +v}dv\eqno (3.38)$$
with $0<\sigma <{1\over 2}-\hbox{\rm Re}\, \alpha -\left|\hbox{\rm Re}
\, \beta\right|$, see [G-R], p. 1015, formula
9.223. For given $y>0$ this function extends to an entire
function of $ $$\left(\alpha ,\beta\right)\in {\bf C}^2$. Indeed, this can be seen from the
formula
$$W_{\alpha ,\beta}\left(y\right)={{y^{\alpha}}\over {e^{y/2}}}\left
({1\over {2\pi i}}\int_{\left(\sigma\right)}{{\Gamma\left(v\right
)\Gamma\left({1\over 2}-\alpha\pm\beta -v\right)}\over {\Gamma\left
({1\over 2}-\alpha\pm\beta\right)}}y^vdv+\sum_{0\le j<-\sigma}{{\left
(-1\right)^j\left({1\over 2}-\alpha\pm\beta\right)_j}\over {j!y^j}}\right
),\eqno (3.39)$$
where $\sigma <{1\over 2}-\hbox{\rm Re}\, \alpha -\left|\hbox{\rm Re}\,
\beta\right|$ and $\sigma$ is not a nonnegative
integer. This is valid for every $\left(\alpha ,\beta\right)\in {\bf C}^
2$ and $y>0$.

We see from (3.39) that if $\alpha$ and $\beta$ are fixed, then for
$0<y<1$ we have $W_{\alpha ,\beta}\left(y\right)\ll_{\delta}y^{\delta}$ for every $
\delta <{1\over 2}-\left|\hbox{\rm Re}\, \beta\right|$.
We also see that $W_{\alpha ,\beta}\left(y\right)$ decays exponentially as $
y\rightarrow\infty$.

The next lemma follows from [G-R], p. 819, 7.625.4 and p. 1022, but since it is very important in our paper we
give a proof of it.

{\bf LEMMA 3.8.} {\it For any} $\hbox{\rm Re$\, S>0$}${\it , for any positive numbers} $
t_1${\it ,} $t_2${\it ,} $M$ {\it and for any}
{\it complex} $k$ {\it and} $\lambda$ {\it such that the function} $
\Gamma$ {\it is regular at the two points} ${1\over 2}-\lambda\pm
it_2${\it , we have that}
$$\int_0^{\infty}y^SW_{k,it_1}(My)W_{\lambda ,it_2}(My){{dy}\over {
y^2}}\eqno (3.40)$$
{\it equals}
$${{M^{1-S}}\over {\Gamma\left({1\over 2}-\lambda\pm it_2\right)}}{
1\over {2\pi i}}\int_{-i\infty}^{i\infty}{{\Gamma\left(-{1\over 2}
\pm it_2+S-s\right)\Gamma\left({1\over 2}\pm it_1+s\right)\Gamma\left
(1-\lambda +s-S\right)}\over {\Gamma\left(1-k+s\right)}}ds,\eqno
(3.41)$$
{\it where the path of integration is chosen in such a way that the poles of the functions}
$\Gamma\left({1\over 2}\pm it_1+s\right)$ {\it and} $\Gamma\left(
1-\lambda +s-S\right)$ {\it lie to the left of the path of integration, and the poles of the functions} $\Gamma\left(-{
1\over 2}\pm it_2+S-s\right)$ {\it lie to the right of it.}

{\it Proof.\/} By a substitution we can assume $M=1$. Using
analytic continuation in $k$, $\lambda$ and $S$ we may assume
$\hbox{\rm Re$\, k=0$}$, Re$\, \lambda =0$, Re$\, S>1$. By these assumptions, using
(3.38) for both Whittaker functions we get that (3.40) equals
$${1\over {2\pi i}}\int_{\left(1/4\right)}{{\Gamma\left(\mu\right
)\Gamma\left({1\over 2}-\lambda\pm it_2-\mu\right)}\over {\Gamma\left
({1\over 2}-\lambda\pm it_2\right)}}I_{\mu}d\mu\eqno (3.42)$$
with
$$I_{\mu}:={1\over {2\pi i}}\int_{\left(1/4\right)}{{\Gamma\left(
v\right)\Gamma\left({1\over 2}-k\pm it_1-v\right)\Gamma\left(S+k+
v+\lambda +\mu -1\right)}\over {\Gamma\left({1\over 2}-k\pm it_1\right
)}}dv,$$
because
$$\int_0^{\infty}e^{-y}y^{S+k+v+\lambda +\mu}{{dy}\over {y^2}}=\Gamma\left
(S+k+v+\lambda +\mu -1\right)$$
by the definition of the $\Gamma$-function. The integral $I_{\mu}$ can
be computed by the $b_4=b_3$ case of the equality of (3.14)
and (3.15), and we get that
$$I_{\mu}={{\Gamma\left(-{1\over 2}\pm it_1+S+\lambda +\mu\right)}\over {
\Gamma\left(S-k+\lambda +\mu\right)}}.$$
Substituting it into (3.42) and writing $s=S+\lambda +\mu -1$ we obtain the lemma.

Part (ii) of the following lemma will be applied directly
in this paper.

{\bf LEMMA 3.9.} {\it (i) Let} $t>0$ {\it be given and let} $n\in
{\bf Z}${\it . We have}
$$\int_0^{\infty}\left|W_{n,it}(y)\Gamma\left({1\over 2}-n\right)\right
|^2{{dy}\over y}\ll_t1.\eqno (3.43)$$
{\it (ii) Let} $t_1,t_2>0$, $S>1/2$ {\it be given and let} $n\in
{\bf Z}$, $M>0$. {\it Then we have that}
$$\int_0^{\infty}\left|y^SW_{n,it_1}(My)W_{0,it_2}(My)\right|{{dy}\over {
y^2}}\ll_{t_1,t_2,S}{{M^{1-S}}\over {\Gamma\left({1\over 2}-n\right
)}}.$$
{\it Proof.\/} To show (i) we apply Lemma 3.8 with $k=\lambda =n$,
$t_1=t_2=t$, $S=M=1$. Note that $W_{n,it}(y)$ is real. Then we
compute (3.41) using the equality of (3.14) and (3.15),
applying it with the parameters
$$a_1={1\over 2}+it,a_2={1\over 2}-it,b_1={1\over 2}-it,b_2=-n,b_
3={1\over 2}+it,b_4=1-n.$$
We get in this way that the left-hand side of (3.43) equals
$$_3F_2\left(\matrix{1,{1\over 2}+it-n,{1\over 2}-it-n\cr
{3\over 2}-it-n,{3\over 2}+it-n\cr}
;1\right){{\Gamma\left({1\over 2}\pm it-n\right)\Gamma\left(1\pm
2it\right)}\over {\Gamma\left({3\over 2}\pm it-n\right)}}.$$
We estimate this series trivially and we get (i).

By a substitution we see that $M=1$ can be assumed in
(ii). The statement then follows from Cauchy-Scwarz,
applying part (i). The lemma is proved.

We finally note that for $x>0$ and arbitrary $\lambda$ and $\mu$ we have
$$x{d\over {dx}}W_{\lambda ,\mu}\left(x\right)=\left(\lambda -{1\over
2}x\right)W_{\lambda ,\mu}\left(x\right)-\left(\mu^2-\left(\lambda
-{1\over 2}\right)^2\right)W_{\lambda -1,\mu}\left(x\right),\eqno
(3.44)$$
see [G-R], p. 1017, 9.234.3.

\noindent{\bf 4. Important lemmas preparing the proof of Theorem 1.1 }
\medskip

In this section $\chi$ will denote a given function satisfying
Condition $D$.

{\bf 4.1 Triple product integrals containing an automorphic kernel function.} Our goal here is to prove Lemma 4.2,
where we give a useful expression for the integral (4.7),
which contains $\overline {B_0(z)}$, a cusp form of weight $0$ and an
automorphic kernel function of weight $1/2$.

{\bf LEMMA 4.1.} {\it Let} $U$ {\it be a cusp form of weight} $0$ {\it for} $
\Gamma_0(4)$. {\it For} $z=x+iy\in\bbb H$ {\it let}
$$V(z):=\overline {B_0(z)}U(z).$$
{\it Let} $k$ {\it be a smooth function on} $[0,\infty )$ {\it such that} $
k^{(j)}(u)$$(u+1)^A$ {\it is bounded on} $[0,\infty )$ {\it for every} $
A>0$ {\it and} $j\ge 0${\it . For} $z,w\in\bbb H$ {\it write}
$$k(z,w):=k\left({{\left|z-w\right|^2}\over {4\hbox{\rm Im$\,z$Im$\,w$}}}\right
)H(z,w)\hbox{\rm \ and }K(z,w):=\sum_{\gamma\in\Gamma_0(4)}k(\gamma
z,w)\overline {\nu (\gamma )}\left({{j_{\gamma}(z)}\over {\left|j_{
\gamma}(z)\right|}}\right)^{-{1\over 2}}.$$
{\it Then for any} $w\in {\bbb H}$ {\it we have }
$$\int_{D_4}V(z)K(z,w)d\mu_z=2\sum_{n=0}^{\infty}\overline {B_n(w
)}\int_{\bbb H}\tilde {k}(Z,i)\overline {\left({{Z-i}\over {Z+i}}\right
)^n}U(T_wZ)d\mu_Z$$
{\it with}
$$\tilde {k}(u):=k(u)\left(u+1\right)^{-{1\over 4}}\qquad (u\in [0
,\infty )),\eqno (4.1)$$
 {\it and}
$$\tilde {k}(z,w):=\tilde {k}\left({{\left|z-w\right|^2}\over {4\hbox{\rm Im$\,
z$Im$\,w$}}}\right)\qquad (z,w\in\bbb H).\eqno (4.2)$$
{\it The sum}
$$\sum_{n=0}^{\infty}\left|\overline {B_n(w)}\int_{\bbb H}\tilde
k(Z,i)\overline {\left({{Z-i}\over {Z+i}}\right)^n}U(T_wZ)d\mu_Z\right
|,$$
{\it as a function of} $w\in D_4${\it , grows at most polynomially at the cusps of} $
\Gamma_0(4)${\it . }

{\it The integral} $\int_{D_4}\overline {V(z)K(z,w)}d\mu_z,$ {\it as a function of} $
w\in D_4${\it , belongs to} $R_{{1\over 2}}(D_4)${\it , and} $\int_{
D_4}\left|V(z)K(z,w)\right|d\mu_z$ {\it decays faster than polynomially at the cusps.}

{\it Proof.\/} It is clear, using (1.3), that if $w\in\bbb H$ is fixed, then for every
$\delta\in\Gamma_0(4)$ and $z\in\bbb H$ we have
$$V(\delta z)=\overline {\nu (\delta )}\left({{j_{\delta}(z)}\over {\left
|j_{\delta}(z)\right|}}\right)^{-{1\over 2}}V(z),\qquad K(\delta
z,w)=\nu (\delta )\left({{j_{\delta}(z)}\over {\left|j_{\delta}(z
)\right|}}\right)^{{1\over 2}}K(z,w).\eqno (4.3)$$
Hence $V(z)K(z,w)$ is invariant in $z$ under $\Gamma_0(4)$, and
$$\int_{D_4}V(z)K(z,w)d\mu_z=\sum_{\gamma\in\Gamma_0(4)}\int_{D_4}
k(\gamma z,w)V(\gamma z)d\mu_z=2\int_{\bbb H}k(z,w)V(z)d\mu_z.\eqno (4.4
)$$
We have $k(T_wZ,T_wi)=k(Z,i)$ by (2.6), because $T_w$ is upper triangular.
Then making the substitution $z=T_wZ$ we get that (4.4) equals
$$2\int_{\bbb H}k(Z,i)V(T_wZ)d\mu_Z.\eqno (4.5)$$
For a given $Z\in\bbb H$ let $Z=i{{1+L}\over {1-L}}$, where $\left
|L\right|<1$. Then it
is easy to see, using also (2.5), that
$$H(Z,i)={{\left(1-L\right)^{{1\over 2}}}\over {\left|1-L\right|^{{
1\over 2}}}},\qquad\hbox{\rm Im}\,\left(T_w\left(i{{1+L}\over {1-L}}\right
)\right)=\left(\hbox{\rm Im$\,w$}\right){{1-\left|L\right|^2}\over {\left
|1-L\right|^2}}.$$
From the definition of $V$ and from (1.2) we then obtain that $k(
Z,i)V(T_wZ)$ equals
$$k\left({{\left|Z-i\right|^2}\over {4\hbox{\rm Im$\,Z$$ $}}}\right
)\left(1-\left|L\right|^2\right)^{{1\over 4}}\overline {\left(\left
(\hbox{\rm Im$\,w$}\right)^{{1\over 4}}\theta (T_wZ)(1-L)^{-{1\over
2}}\right)}U(T_wZ).$$
It is easy to check that
$$\qquad\left(1-\left|L\right|^2\right)^{{1\over 4}}=\left(1+{{\left
|Z-i\right|^2}\over {4\hbox{\rm Im$\,Z$$ $}}}\right)^{-{1\over 4}},
\qquad L={{Z-i}\over {Z+i}}.$$
So, taking the Taylor expansion (2.2) for $w$ in place of $z$,
we get for every $Z\in\bbb H$ that
$$k(Z,i)V(T_wZ)=\tilde {k}(Z,i)U(T_wZ)\sum_{n=0}^{\infty}\overline {
B_n(w)}\overline {\left({{Z-i}\over {Z+i}}\right)^n}.\eqno (4.6)$$
We need the weak estimate that if $w\in D_4$ and $0\le j\le 5$ is
such that $\gamma_j^{-1}w\in D_1$, then $\left|B_n(w)\right|\ll\left
(n+1\right)^{A_0}\left(\hbox{\rm Im}\,\gamma_j^{-1}w\right)^{A_0}$
with some absolute constant $A_0$. This follows easily from
Lemma 6.2 of [B3]. Using that $U$ is bounded, we then see
that inserting (4.6) into (4.5) we can integrate term by
term. In this way we get the assertions of the lemma except the last sentence. In the last sentence the transformation
property follows easily from (2.6). For the estimates we use
$$\left|K(z,w)\right|\le\sum_{\gamma\in SL(2,{\bf Z})}\left|k\left
({{\left|\gamma z-w\right|^2}\over {4\hbox{\rm Im$\,\gamma z$Im$\,w$}}}\right
)\right|,$$
and we note that if $z,w\in D_1$, then for any $A>0$ the right-hand side
here
$$\ll_A\left(\hbox{\rm Im}\,z\right)^{A_0}\left(1+{{{\rm I}{\rm m}\,w}\over {
{\rm I}{\rm m}\,z}}\right)^{-A}$$
with some absolute constant $A_0$. This follows from Lemma 6.3 of [B3]. Using that $
V$
decays faster than polynomially at the cusps, the lemma
follows.

{\bf LEMMA 4.2.} {\it Let} $k=k_{\chi}$ {\it be the function defined in Lemma 3.7, and let} $
K(z,w)$ {\it be as in Lemma 4.1. Let} $u$ {\it be a cusp form of weight} $
0$ {\it for} $SL(2,{\bf Z})$
{\it with} $\Delta_0u=S\left(S-1\right)u$, $S={1\over 2}+iT${\it . Then for any} $
w\in\bbb H$ {\it we have that}
$$\int_{D_4}\overline {B_0(z)}u(4z)K(z,w)d\mu_z\eqno (4.7)$$
{\it equals}
$${4\over {\Gamma\left({1\over 2}\pm iT\right)}}\sum_{n=0}^{\infty}{{
C_{n,\chi}\left(T\right)}\over {\Gamma\left({1\over 2}+n\right)}}\overline {
B_n\left(w\right)}\left(K_{n-1}K_{n-2}\ldots K_1K_0u\right)(4w),\eqno
(4.8)$$
{\it where}
$$C_{n,\chi}\left(T\right):=\int_{-\infty}^{\infty}\left|{{\Gamma\left
({1\over 4}+it\right)\Gamma\left({1\over 4}+it\pm iT\right)}\over {
\Gamma\left(2it\right)}}\right|^2{{S_n\left(t^2;{1\over 4}+iT,{1\over
4},{1\over 4}-iT\right)}\over {\left({1\over 2}+iT\right)_n\left({
1\over 2}-iT\right)_n}}\chi (t)dt.$$
{\it The sum} $\sum_{n=0}^{\infty}\left|{{C_{n,\chi}\left(T\right
)}\over {\Gamma\left({1\over 2}+n\right)}}\overline {B_n\left(w\right
)}\left(K_{n-1}K_{n-2}\ldots K_1K_0u\right)(4w)\right|,$ {\it as a function of} $
w\in D_4${\it , grows at most polynomially at the cusps.}

{\it Proof.\/} We will apply Lemma 4.1 with $U(z)=u(4z)$, and we
use the notations (4.1), (4.2). Remark that if $n\ge 0$, $w\in\bbb
H$,
then we get
$$\int_{\bbb H}\tilde {k}(Z,i)\overline {\left({{Z-i}\over {Z+i}}\right
)^n}u(4T_wZ)d\mu_Z=\int_{\bbb H}\tilde {k}(z,w)\overline {\left({{
z-w}\over {z-\overline w}}\right)^n}u(4z)d\mu_z\eqno (4.9)$$
by the substitution $z=T_wZ$. We now make a transition to geodesic polar
coordinates around $w$, i.e. we use (2.1) with $w$ in place of
$z_0$. See also the form of the invariant measure given below (2.1). We get in this way that (4.9) equals
$$\int_0^{\infty}\tilde {k}\left({{\tanh^2({r\over 2})}\over {1-\tanh^
2({r\over 2})}}\right)\tanh^n({r\over 2}^{})\left(\int_0^{2\pi}u(
4z)e^{-in\phi}d\phi\right)\sinh rdr.\eqno (4.10)$$
Here we write $R$ in place of $\tanh({r\over 2})$, and use
$$\sinh rdr={{4R}\over {\left(1-R^2\right)^2}}dR.\eqno (4.11)$$
We apply also Lemma 2.2 and (4.1), and we get in this
way that (4.10) is the same as
$${{8\pi I_k\left(S\right)}\over {n!}}\left(K_{n-1}K_{n-2}\ldots
K_1K_0u\right)(4w),$$
 where
$$I_k\left(S\right):=\int_0^1\hbox{\rm $k\left({{R^2}\over {1-R^2}}\right
)R^{2n+1}\left(1-R^2\right)^{-{7\over 4}}F\left(S,1-S,n+1;{{R^2}\over {
R^2-1}}\right)dR$}.$$
Using here the definition of $k=k_{\chi}$ from Lemma 3.7 the resulting double
integral is easily seen to be absolutely convergent. We
get, writing $u={{R^2}\over {1-R^2}}$ and applying (3.5) that
$$I_k\left(S\right)={1\over {2\pi}}\int_0^{\infty}\left|{{\Gamma\left
({1\over 4}+it\right)\Gamma\left({3\over 4}+it\right)}\over {\Gamma\left
(2it\right)}}\right|^2\chi (t)I_n\left(t,T\right)dt,$$
where
$$I_n\left(t,T\right):=\int_0^{\infty}F\left({3\over 4}-it,{3\over
4}+it,1;-u\right)u^nF({1\over 2}+iT+n,{1\over 2}-iT+n,1+n;-u)du.$$
Using Lemma 4.1, Lemma 3.4 and (3.6) we get that (4.7) equals (4.8). The last statement of the lemma follows from the
corresponding statement of Lemma 4.1.

{\bf 4.2. An expression for the spectral sum.} We give an
expression for the spectral sum of Theorem 1.1 in terms
of an automorphic kernel function. Our main result here
is Lemma 4.5.

{\bf LEMMA 4.3.} {\it The notations and assumptions of Lemma 4.1 are valid. Let} $
f$ {\it be a Maass form of weight} ${1\over 2}$ {\it for} $\Gamma_
0(4)$ {\it with }
$\Delta_{1/2}f=s(s-1)f$ {\it for some\/} Re$\, s\ge{1\over 2}$, $s={
1\over 2}+it$. {\it Then}
$$\int_{D_4}\left(\int_{D_4}V(z)K(z,w)d\mu_z\right)f(w)d\mu_w=16\pi\left
(\int_{D_4}V(z)f(z)d\mu_z\right)L_k(s),\eqno (4.12)$$
 {\it where}
$$L_k(s):=\int_0^1k\left({{R^2}\over {1-R^2}}\right)\left(1-R^2\right
)^{-9/4}F\left(s+{1\over 4},{5\over 4}-s,1;{{R^2}\over {R^2-1}}\right
)RdR.$$
{\it If} $k=k_{\chi}$ {\it is the function defined in Lemma 3.7, then}
$$L_k(s)=\chi\left(t\right).\eqno (4.13)$$
{\it Proof.\/} Taking real and imaginary parts, we may assume that
$k(u)$ is real for any $u\in [0,\infty )$. Since $k$ is real, it
is not hard to see, using (2.6), (2.7) and (1.3) that
$K(z,w)=\overline {K(w,z)}$. Hence by (4.3) and the transformation formulas satisfied by $
f$ we
see that $K(z,w)f(w)$ is invariant in $w$ under $\Gamma_0(4)$. We
also see that
$$\int_{D_4}\left(\int_{D_4}V(z)K(z,w)d\mu_z\right)f(w)d\mu_w=\int_{
D_4}V(z)\left(\int_{D_4}f(w)\overline {K(w,z)}d\mu_w\right)d\mu_z
,$$
the application of the Fubini theorem is justified by
the last statement of Lemma 4.1. By the definition of $K$
we see that
$$\int_{D_4}f(w)\overline {K(w,z)}d\mu_w=2\int_{\bbb H}f(w)\overline {
k(w,z)}d\mu_w.\eqno (4.14)$$
We have
$$\overline {H(w,z)}=i^{{{-1}\over 2}}\left({{w-\overline z}\over {\left
|w-\overline z\right|}}\right)^{{1\over 2}}=\left({{w-\overline z}\over {
z-\overline w}}\right)^{{1\over 4}},$$
the last equality holds because the fourth powers are the same, and the arguments of
both sides lie in $(-{{\pi}\over 4},{{\pi}\over 4})$. We use geodesic polar coordinates
around $z$ (see (2.1)) and we write
$$F(r,\phi ):=f(w)\left({{w-\overline z}\over {z-\overline w}}\right
)^{{1\over 4}}.$$
We get in this way that (4.14) equals
$$2\int_0^{\infty}k\left({{\tanh^2({r\over 2})}\over {1-\tanh^2({
r\over 2})}}\right)\left(\int_0^{2\pi}F(r,\phi )d\phi\right)\sinh
rdr.$$
By Lemma 2.2 we get
$$\int_0^{2\pi}F(r,\phi )d\phi =2\pi f(z)\left(1-\tanh^2({r\over
2})\right)^{-1/4}F\left(s+{1\over 4},{5\over 4}-s,1;{{\tanh^2({r\over
2})}\over {\tanh^2({r\over 2})-1}}\right).$$
Writing $R$ in place of $\tanh({r\over 2})$, using (4.11) we obtain (4.12). By the substitution $
u={{R^2}\over {1-R^2}}$ and by Lemma 3.7 we get (4.13), the lemma is proved.

{\bf LEMMA 4.4.} {\it If} $f_1,f_2\in R_{{1\over 2}}(D_4)${\it , then we have that} $\left
(f_1,f_2\right)_4$
{\it equals}
$$\sum_{j=0}^{\infty}\left(f_1,u_{j,{1\over 2}}\right)_4\overline {\left
(f_2,u_{j,{1\over 2}}\right)_4}+{1\over {4\pi}}\sum_{\eufm{a}=0,\infty}
\int_{-\infty}^{\infty}\zeta_{\eufm{a}}(f_1,r)\overline {\zeta_{\eufm{
a}}(f_2,r)}dr.$$
{\it Proof.\/} This is well-known, see [P], formula (27).

{\bf LEMMA 4.5.} {\it Let} $k=k_{\chi}$ {\it be the function defined in Lemma 3.7, and let} $
K(z,w)$ {\it be as in Lemma 4.1. Let} $u_1$ {\it and} $u_2$
{\it be two cusp forms of weight} $0$ {\it for} $SL(2,{\bf Z})${\it . Then }
$$\int_{D_4}\left(\int_{D_4}\overline {B_0(z)}u_1\left(4z\right)K
(z,w)d\mu_z\right)B_0\left(w\right)\overline {u_2\left(4w\right)}
d\mu_w\eqno (4.15)$$
 {\it equals the sum of}
$$16\pi\sum_{j=1}^{\infty}\chi\left(T_j\right)\left(B_0\kappa\left
(\overline {u_2}\right),u_{j,{1\over 2}}\right)_4\overline {\left
(B_0\kappa\left(\overline {u_1}\right),u_{j,{1\over 2}}\right)_4}\eqno
(4.16)$$
 {\it and}
$$4\sum_{\eufm{a}=0,\infty}\int_{-\infty}^{\infty}\chi\left(r\right
)\zeta_{\eufm{a}}\left(B_0\kappa\left(\overline {u_2}\right),r\right
)\overline {\zeta_{\eufm{a}}\left(B_0\kappa\left(\overline {u_1}\right
),r\right)}dr.\eqno (4.17)$$
 {\it Proof.\/} Let
$$f_1\left(w\right):=B_0\left(w\right)\overline {u_2\left(4w\right
)},\quad f_2\left(w\right):=\int_{D_4}B_0(z)\overline {u_1\left(4z\right
)K(z,w)}d\mu_z.$$
If $f$ is a Maass form of weight ${1\over 2}$ for $\Gamma_0(4)$ with
$\Delta_{1/2}f=s(s-1)f$ for some $s={1\over 2}+it$, then we have by
Lemma 4.3 that
$$\overline {\left(f_2,f\right)_4}=16\pi\chi\left(t\right)\int_{D_
4}\overline {B_0(z)}u_1\left(4z\right)f(z)d\mu_z.$$
Lemma 4.4 implies that (4.15) equals the sum of (4.16) and
(4.17), but at the moment it seems that $j=0$ should be present in
the summation in (4.16). However, that term is $0$ by Lemma 6.6 of [B3]. The lemma is proved.

{\bf 4.3. Writing} $B_t\overline {B_0}$  {\bf as an Eisenstein series.} We
mentioned in Sections 1.4 that it is very important for
our proof that the functions $B_t\overline {B_0}$ are linear combinations
of Eisenstein series. We write a certain average of this function as an incomplete Eisenstein
series in Lemma 4.7.

{\bf LEMMA 4.6.}  {\it For} $
z\in\bbb H$ {\it let}
$$F(z):=\sum_{j=0}^5\left|B_0(\gamma_jz)\right|^2,\qquad G(z):=\sum_{
\gamma\in\Gamma_{\infty}\setminus SL(2,{\bf Z})}\psi\left(\hbox{\rm Im$\,
(\gamma z)$}\right),$$
{\it where }
$$\psi\left(y\right):=\sum_{m=1}^{\infty}e^{-\pi{{m^2}\over y}}.$$
 {\it Then for every} $z\in\bbb H$ {\it we have} $F(z)=6G(z)
+3.$

{\it Proof.\/} During the proof of Lemma 6.6 of [B3] (see the
last lines of p. 632) it is shown that $F(z)=DG(z)+C$ for $z\in\bbb
H$ with some constants $C$
and $D$. So it is enough to determine these constants.

Recall the definitions of $\gamma_j$ from Section 2.5. Note first
that $B_0\left(Z-{1\over 2}\right)=\sqrt 2B_0\left(4Z\right)-B_0\left
(Z\right)$ for $Z\in\bbb H$ by (1.2). One has $\gamma_5z=-{1\over
w}-{1\over 2}$ with
$w=4z-2$, hence using also (1.4) we get
$$B_0\left(\gamma_5z\right)=\sqrt 2B_0\left(-{4\over w}\right)-B_
0\left(-{1\over w}\right)=e\left({{-1}\over 8}\right)\left({w\over {\left
|w\right|}}\right)^{{1\over 2}}\left(\sqrt 2B_0\left({w\over {16}}\right
)-B_0\left({w\over 4}\right)\right)$$
for every $z\in\bbb H$. This shows by (1.2) that
$B_0\left(\gamma_5\left(iy\right)\right)=o\left(1\right)$ as $y\rightarrow
\infty$. For $0\le j\le 3$ it is clear by
(1.4) that we have $\left|B_0(\gamma_jz)\right|=\left|B_0\left({{
z+j}\over 4}\right)\right|$. We easily get
from these remarks and (1.2) that
$F(iy)=3y^{1/2}+o\left(1\right)$ as $y\rightarrow\infty$. On the other hand, it is easy to see that
$G(iy)=\psi\left(y\right)+o\left(1\right)$ as $y\rightarrow\infty$, and it follows from (1.4) that
$$1+2\psi\left(y\right)=y^{1/2}\sum_{m=-\infty}^{\infty}e^{-\pi m^
2y}=y^{1/2}+o\left(1\right).$$
Letting $y\rightarrow\infty$ we get the lemma.

{\bf LEMMA 4.7.} {\it If} $t\ge 0$ {\it is an integer, for} $z\in\bbb
H$ {\it let}
$$F_t(z):=\sum_{j=0}^5B_t(\gamma_jz)\overline {B_0(\gamma_jz)}\left
({{j_{\gamma_j}\left(z\right)}\over {\left|j_{\gamma_j}\left(z\right
)\right|}}\right)^{-2t}$$
{\it and }
$$G_t(z):=\sum_{\gamma\in\Gamma_{\infty}\setminus SL(2,{\bf Z})}\psi_
t\left(\hbox{\rm Im$\,(\gamma z)$}\right)\left({{j_{\gamma}\left(z\right
)}\over {\left|j_{\gamma}\left(z\right)\right|}}\right)^{-2t},$$
{\it where}
$$\psi_t\left(y\right):={1\over {t!}}\sum_{m=1}^{\infty}e^{-\pi{{m^
2}\over y}}\left({{\pi m^2}\over y}\right)^t.$$
{\it Then for every} $t>0$ {\it we have} $F_t(z)=6G_t(z).$

{\it Proof.\/} We have $F_0=F$, $G_0=G$ (see Lemma 4.6). It is easy to see that for
every $t\ge 0$ we have
$${1\over {t+1}}K_tG_t=G_{t+1},$$
this follows from the identity
$${1\over {t+1}}\left(\psi_t^{(1)}\left(y\right)y+t\psi_t\left(y\right
)\right)=\psi_{t+1}\left(y\right)$$
and Lemma 2.1 (vi). Using Lemma 4.6 we see that it is
enough to prove that
$${1\over {t+1}}K_tF_t=F_{t+1}\eqno (4.18)$$
for every $t\ge 0$. We use Lemma 2.1 (i) with $k_1:=t+{1\over 4}$, $
k_2:={1\over 4}$,
$$f(z):=B_t(\gamma_jz)\left({{j_{\gamma_j}\left(z\right)}\over {\left
|j_{\gamma_j}\left(z\right)\right|}}\right)^{-2t-{1\over 2}},\qquad
g(z):=B_0(\gamma_jz)\left({{j_{\gamma_j}\left(z\right)}\over {\left
|j_{\gamma_j}\left(z\right)\right|}}\right)^{-{1\over 2}}.$$
Then $K_{-k_2}\left(\overline g\right)=0$ by Lemma 2.1 (vi) and (v). So (4.18) follows using (2.3) and Lemma 2.1 (vi).
The lemma is proved.

\noindent{\bf 5. Proof of the theorem}
\medskip

{\bf 5.1. A special case.} We first assume that $\chi$ is a function satisfying Condition $
D$.

By Lemma 4.5 and Lemma 4.2 we have that the sum of
$$16\pi\sum_{j=1}^{\infty}\chi\left(T_j\right)\left(B_0\kappa\left
(\overline {u_2}\right),u_{j,{1\over 2}}\right)_4\overline {\left
(B_0\kappa\left(\overline {u_1}\right),u_{j,{1\over 2}}\right)_4}\eqno
(5.1)$$
 and
$$4\sum_{\eufm{a}=0,\infty}\int_{-\infty}^{\infty}\chi\left(r\right
)\zeta_{\eufm{a}}\left(B_0\kappa\left(\overline {u_2}\right),r\right
)\overline {\zeta_{\eufm{a}}\left(B_0\kappa\left(\overline {u_1}\right
),r\right)}dr\eqno (5.2)$$
 equals
$${4\over {\Gamma\left({1\over 2}\pm it_1\right)}}\sum_{n=0}^{\infty}{{
C_{n,\chi}\left(t_1\right)}\over {\Gamma\left({1\over 2}+n\right)}}
J_n(u_1,u_2),\eqno (5.3)$$
 where
$$J_n=J_n(u_1,u_2):=\int_{D_4}B_0(w)\overline {B_n(w)}\overline {u_
2\left(4w\right)}\left(K_{n-1}K_{n-2}\ldots K_1K_0u_1\right)(4w)d
\mu_w.$$
Let us write
$$f\left(z\right)=f_{n,u_1,u_2}\left(z\right):=\overline {u_2\left
(z\right)}\left(K_{n-1}\ldots K_1K_0u_1\right)\left(z\right).\eqno
(5.4)$$
We then have that
$$f(\gamma z)=\left({{j_{\gamma}(z)}\over {\left|j_{\gamma}(z)\right
|}}\right)^{2n}f(z)\eqno (5.5)$$
for every $\gamma\in SL(2,{\bf Z})$. Since the substitution $w\rightarrow
-{1\over {4w}}$
normalizes $\Gamma_0(4)$, so
$$J_n=\int_{D_4}B_0\left(w\right)\overline {B_n(w)}f\left(4w\right
)d\mu_w=\int_{D_4}B_0\left({{-1}\over {4w}}\right)\overline {B_n\left
({{-1}\over {4w}}\right)}f\left({{-1}\over w}\right)d\mu_w,\eqno
(5.6)$$
hence by (2.4) and (5.5) we get
$$J_n=\int_{D_4}B_0\left(w\right)\overline {B_n(w)}f\left(w\right
)d\mu_w.$$
Using again (5.5), we finally get
$$J_n=\int_{D_1}\overline {F_n(w)}f\left(w\right)d\mu_w,\eqno (5.
7)$$
with the function $F_n$ defined in Lemma 4.7. Using Lemmas
4.6 and 4.7 we see by unfolding that
$$J_n=6\int_0^{\infty}\int_0^1\overline {\psi_n(y)}f\left(x+iy\right
){{dxdy}\over {y^2}}+3\delta_{0,n}\int_{D_1}f\left(w\right)d\mu_w
,\eqno (5.8)$$
where $\delta_{0,n}$ is Kronecker's symbol.

It is trivial by our assumptions that if $n=0$, then
$$\int_{D_1}f\left(w\right)d\mu_w=\delta_{u_1,u_2}.\eqno (5.9)$$
It is well-known that if $u$ is a cusp form of weight $0$ for $SL
(2,{\bf Z})$ with $\Delta_0u$$=s(s-1)u$, where
$s={1\over 2}+it$, and
$$u(z)=\sum_{m\neq 0}\rho_u(m)W_{0,it}(4\pi\left|m\right|y)e(mx),$$
then for any $N\ge 0$ we have
$$\left(K_{N-1}K_{N-2}\ldots K_1K_0u\right)(z)=\sum_{m\neq 0}\rho_
N^u(m)W_{N{\rm s}{\rm g}{\rm n}\left(m\right),it}(4\pi\left|m\right
|y)e(mx)$$
with
$$\rho_N^u(m)=(-1)^N\rho_u(m)\eqno (5.10)$$
for $m>0$, and
$$\rho_N^u(m)=\left(s\right)_N\left(1-s\right)_N\rho_u(m)\eqno (5
.11)$$
for $m<0$, but we show now these statements for the
sake of completeness. Indeed, by (3.44), if $m>0$, then
$$L_k\left(W_{k,it}(4\pi my)e(mx)\right)=-\left(\left(it\right)^2
-\left(k-{1\over 2}\right)^2\right)W_{k-1,it}(4\pi my)e(mx),\eqno
(5.12)$$
and if $m<0$, then
$$K_k\left(W_{-k,it}(4\pi\left|m\right|y)e(mx)\right)=-\left(\left
(it\right)^2-\left(-k-{1\over 2}\right)^2\right)W_{-k-1,it}(4\pi\left
|m\right|y)e(mx).\eqno (5.13)$$
Since, by Lemma 2.1 (iv), we have
$$\left(L_1L_2\ldots L_N\left(K_{N-1}K_{N-2}\ldots K_1K_0u\right)\right
)(z)={{\Gamma (s+N)}\over {\Gamma (s-N)}}u(z),$$
so by repeated application of (5.12) we get (5.10), and by repeated application
of (5.13) we get (5.11).

It is easy to see using (5.4) that
$$\int_0^{\infty}\int_0^1\overline {\psi_n(y)}f\left(x+iy\right){{
dxdy}\over {y^2}}=\sum_{m\neq 0}\overline {\rho_{u_2}(m)}\rho_n^{
u_1}(m)I_{n,t_1,t_2}\left(m\right)\eqno (5.14)$$
with
$$I_{n,t_1,t_2}\left(m\right):=\int_0^{\infty}\psi_n(y)W_{n{\rm s}{\rm g}{\rm n}\left(m\right),it_1}
(4\pi\left|m\right|y)W_{0,it_2}(4\pi\left|m\right|y){{
dy}\over {y^2}}\eqno (5.15)$$
(remark that $\psi_n(y)$ is real). By the well-known formula
$${1\over {2\pi i}}\int_{(\sigma )}\Gamma\left(S\right)Y^{-S}dS=e^{
-Y}$$
we see for every $l\ge 0$ and $\sigma >{1\over 2}$ that
$$\psi_l(y)={1\over {l!}}{1\over {2\pi i}}\int_{(\sigma )}\pi^{-S}
\zeta\left(2S\right)\Gamma\left(l+S\right)y^SdS.\eqno (5.16)$$
We will compute (5.3) by (5.8), (5.14), (5.15), (5.16),
in this way we get summations over $m,n$ and integration
over $y$ and $S$. We can see that if $\sigma$ is fixed to be a
large enough absolute constant, then these summations and
integrations are absolutely convergent together. This can
be seen by the definition of $C_{n,\chi}$ in Lemma 4.2, by Lemma 3.3, (5.10), (5.11), estimating the
integral involving Whittaker functions by Lemma 3.9 (ii).

Applying Lemma 3.8, we get for any
$\hbox{\rm Re$\, S>0$}$ that
$$\int_0^{\infty}y^SW_{n{\rm s}{\rm g}{\rm n}\left(m\right),it_1}
(4\pi\left|m\right|y)W_{0,it_2}(4\pi\left|m\right|y){{dy}\over {y^
2}}\eqno (5.17)$$
equals
$${{\left(4\pi\left|m\right|\right)^{1-S}}\over {\Gamma\left({1\over
2}\pm it_2\right)}}{1\over {2\pi i}}\int_{-i\infty}^{i\infty}{{\Gamma\left
(-{1\over 2}\pm it_2+S-s\right)\Gamma\left({1\over 2}\pm it_1+s\right
)\Gamma\left(1+s-S\right)}\over {\Gamma\left(1-n+s\right)}}ds\eqno
(5.18)$$
 in the case $m>0$, and
$${{\left(4\pi\left|m\right|\right)^{1-S}}\over {\Gamma\left({1\over
2}+n\pm it_1\right)}}{1\over {2\pi i}}\int_{-i\infty}^{i\infty}{{
\Gamma\left(-{1\over 2}\pm it_1+S-s\right)\Gamma\left({1\over 2}\pm
it_2+s\right)\Gamma\left(1+n+s-S\right)}\over {\Gamma\left(1+s\right
)}}ds\eqno (5.19)$$
in the case $m<0$. Indeed, we obtain it by the choice
$$k=n,\lambda =0\hbox{\rm \ in the case $m>0,$}$$
$$k=0,\lambda =-n\hbox{\rm \ in the case $m<0.$}$$
In the case $m<0$ we apply Lemma 3.8 by exchanging $t_1$
and $t_2$.

By (5.8), (5.9) and (5.14) we have
$$J_n=3\delta_{0,n}\delta_{u_1,u_2}+6\sum_{m\neq 0}\overline {\rho_{
u_2}(m)}\rho_n^{u_1}(m)I_{n,t_1,t_2}\left(m\right).\eqno (5.20)$$
We can determine $\rho_n^{u_1}(m)$ by (5.10) and (5.11). We use
that
$$\rho_{u_1}(m)\overline {\rho_{u_2}(m)}=\rho_{u_1}(-m)\overline {
\rho_{u_2}(-m)}$$
for every $m\neq 0$, since it is assumed that either $u_1$ and $u_
2$ are odd, or both of them are even. We then
see by (5.15), (5.16), (5.17), (5.18), (5.19) and (1.6) that fixing $
\sigma$ to be a large enough absolute constant,
$$\Gamma\left({1\over 2}\pm it_1\right)\sum_{n=0}^{\infty}{{C_{n,
\chi}\left(t_1\right)}\over {\Gamma\left({1\over 2}+n\right)}}\sum_{
m\neq 0}\overline {\rho_{u_2}(m)}\rho_n^{u_1}(m)I_{n,t_1,t_2}\left
(m\right)\eqno (5.21)$$
equals
$$4\pi{1\over {2\pi i}}\int_{(\sigma )}\left(4\pi^2\right)^{-S}\zeta\left
(2S\right)L\left(S\right)\left(E^{+}\left(S,{1\over 2}\right)+{{\sin
\pi s_2}\over {\sin\pi s_1}}E^{-}\left(S,{1\over 2}\right)\right)
dS,\eqno (5.22)$$
where $E^{+}\left(S,D\right)$ denotes the sum
$$\sum_{n=0}^{\infty}{{C_{n,\chi}\left(t_1\right)\Gamma\left(S+n\right
)}\over {n!\Gamma\left(D+n\right)2\pi i}}\int_{(\tau )}{{\Gamma\left
(-{1\over 2}\pm it_1+S-s\right)\Gamma\left({1\over 2}\pm it_2+s\right
)\Gamma\left(1+n+s-S\right)}\over {\Gamma\left(1+s\right)}}ds,$$
and $E^{-}\left(S,D\right)$ denotes the sum
$$\sum_{n=0}^{\infty}(-1)^n{{C_{n,\chi}\left(t_1\right)\Gamma\left
(S+n\right)}\over {n!\Gamma\left(D+n\right)2\pi i}}\int_{(\tau )}{{
\Gamma\left(-{1\over 2}\pm it_2+S-s\right)\Gamma\left({1\over 2}\pm
it_1+s\right)\Gamma\left(1+s-S\right)}\over {\Gamma\left(1-n+s\right
)}}ds\eqno (5.23)$$
with $\hbox{\rm Re$\, S$}-1<\tau <\hbox{\rm Re$\, S$}-{1\over 2}$, $\tau
>-{1\over 2}$. There is such a$ $ $\tau$ for every
$\hbox{\rm Re$\, S$}>0$. Our computations are justified by the discussion
below (5.16). We see by Lemma 3.3 that the summation and
integrations in $n,s$ and $S$ are absolutely convergent.

See Section 2.4 for the properties of the function
$\zeta\left(2S\right)L\left(S\right)$. This function is regular at $
S=1$ if $u_1\neq u_2$, and its
residue at $S=1$ in the case $u_1=u_2$ is
$$\hbox{\rm res}_{S=1}\zeta\left(2S\right)L\left(S\right)=\zeta\left
(2\right){1\over {\Gamma\left({1\over 2}\pm it_1\right)\left|D_1\right
|}}={{\pi}\over {2\Gamma\left({1\over 2}\pm it_1\right)}}.\eqno (
5.24)$$
This follows from [I], (8.12), (8.9) and (8.5), taking into account $\zeta\left(2\right)={{\pi^2}\over 6},\;\left
|D_1\right|={{\pi}\over 3}.$ The last relation follows from [I], (6.33), (3.26).

Since $C_{n,\chi}\left(t_1\right)$ decreases faster than polynomially in $
n$
by Lemma 3.3, so by the properties of $L\left(S\right)$ we see shifting the
integration to the left that for a small $\epsilon >0$, e.g. take
$\epsilon ={1\over {100}}$, we have that (5.22) equals the sum of
$$4\pi{1\over {2\pi i}}\int_{\left({1\over 2}-\epsilon\right)}\left
(4\pi^2\right)^{-S}\zeta\left(2S\right)L\left(S\right)\left(E^{+}\left
(S,{1\over 2}\right)+{{\sin\pi s_2}\over {\sin\pi s_1}}E^{-}\left
(S,{1\over 2}\right)\right)dS\eqno (5.25)$$
and
$$\delta_{u_1,u_2}{1\over {2\Gamma\left({1\over 2}\pm it_1\right)}}\left
(E^{+}\left(1,{1\over 2}\right)+E^{-}\left(1,{1\over 2}\right)\right
).\eqno (5.26)$$
Note that this last term is present only in the case
$t_1=t_2$. We now determine $E^{+}\left(1,{1\over 2}\right)+E^{-}\left
(1,{1\over 2}\right)$ in the case
$t_1=t_2$. Since $(-1)^n{{\Gamma\left(s\right)}\over {\Gamma\left
(1-n+s\right)}}=-{{\Gamma\left(n-s\right)}\over {\Gamma\left(1-s\right
)}},$ so, using the substitution $s\rightarrow -s$ in the integral in
(5.23), we have that $E^{+}\left(1,{1\over 2}\right)+E^{-}\left(1
,{1\over 2}\right)$ equals
$$\sum_{n=0}^{\infty}{{C_{n,\chi}\left(t_1\right)}\over {\Gamma\left
({1\over 2}+n\right)2\pi i}}\left(\int_{(\tau )}{{\Gamma\left({1\over
2}\pm it_1\pm s\right)\Gamma\left(n+s\right)}\over {\Gamma\left(1
+s\right)}}ds-\int_{(-\tau )}{{\Gamma\left({1\over 2}\pm it_1\pm
s\right)\Gamma\left(n+s\right)}\over {\Gamma\left(1+s\right)}}ds\right
)$$
with $0<\tau <{1\over 2}$. For $n>0$ the difference of these
integrals is 0, and for $n=0$ it is $2\pi i$$\Gamma^2\left({1\over
2}\pm it_1\right)$. Hence,
if $t_1=t_2$, we have
$$E^{+}\left(1,{1\over 2}\right)+E^{-}\left(1,{1\over 2}\right)={{
C_{0,\chi}\left(t_1\right)}\over {\Gamma\left({1\over 2}\right)}}
\Gamma^2\left({1\over 2}\pm it_1\right).\eqno (5.27)$$
It is clear that (5.25) equals
$$\lim_{\delta\rightarrow 0+0}{{4\pi}\over {2\pi i}}\int_{\left({
1\over 2}-\epsilon\right)}e^{\delta S^2}\left(4\pi^2\right)^{-S}\zeta\left
(2S\right)L\left(S\right)\left(E^{+}\left(S,{1\over 2}\right)+{{\sin
\pi s_2}\over {\sin\pi s_1}}E^{-}\left(S,{1\over 2}\right)\right)
dS,\eqno (5.28)$$
and for a given $\delta >0$ we have that
$${{4\pi}\over {2\pi i}}\int_{\left({1\over 2}-\epsilon\right)}e^{
\delta S^2}\left(4\pi^2\right)^{-S}\zeta\left(2S\right)L\left(S\right
)\left(E^{+}\left(S,D\right)+{{\sin\pi s_2}\over {\sin\pi s_1}}E^{
-}\left(S,D\right)\right)dS\eqno (5.29)$$
is a regular function of $D$ for $\hbox{\rm Re$\, D>0$}$. Let us consider this
function first for large enough Re$\, D$. By the definition of
$E^{+}\left(S,D\right)$, $E^{-}\left(S,D\right)$, $C_{n,\chi}\left
(t_1\right)$, the upper bound for $\chi$ and
Lemma 3.1 we see that if $D$ has large enough real part, then we can compute
(5.29) by inserting the defining integral for
$C_{n,\chi}\left(t_1\right)$ in $E^{+}\left(S,D\right)$ and $E^{-}\left
(S,D\right)$, since the resulting triple integral in $s,S,t$ and
summation in $n$ are absolutely convergent. We will use
the following two identities, both of them follow from
Lemma 3.2.

For any $S$ and $s$ with $\hbox{\rm Re$\, S$}={1\over 2}-\epsilon$, $\hbox{\rm Re$\,
s$}=-{1\over 4}-{{\epsilon}\over 2}$ we have
that
$$\sum_{n=0}^{\infty}{{_3F_2\left(\matrix{-n,{1\over 4}+it,{1\over
4}-it\cr
{1\over 2}+it_1,{1\over 2}-it_1\cr}
;1\right)}\over {n!}}{{\Gamma\left(S+n\right)\Gamma\left(1-S+s+n\right
)}\over {\Gamma\left(D+n\right)}}$$
equals
$${{\Gamma\left({1\over 2}\pm it_1\right)}\over {\Gamma\left({1\over
4}\pm it\right)\Gamma\left(D-S\right)\Gamma\left(D+S-s-1\right)}}
F_1\left(s\right),$$
defining $F_1\left(s\right)$ as
$${1\over {2\pi i}}\int_{\left(-c\right)}{{\Gamma\left({1\over 4}
\pm it+T\right)\Gamma\left(S+T\right)\Gamma\left(1-S+s+T\right)\Gamma\left
(-T\right)\Gamma\left(D-1-s-T\right)}\over {\Gamma\left({1\over 2}
\pm it_1+T\right)}}dT$$
with ${1\over 4}-{{\epsilon}\over 2}<c<{1\over 4}$; and
$$\sum_{n=0}^{\infty}{{_3F_2\left(\matrix{-n,{1\over 4}+it,{1\over
4}-it\cr
{1\over 2}+it_1,{1\over 2}-it_1\cr}
;1\right)}\over {n!}}{{\Gamma\left(S+n\right)\Gamma\left(-s+n\right
)}\over {\Gamma\left(D+n\right)}}$$
equals
$${{\Gamma\left({1\over 2}\pm it_1\right)}\over {\Gamma\left({1\over
4}\pm it\right)\Gamma\left(D-S\right)\Gamma\left(D+s\right)}}F_2\left
(s\right),$$
defining, again with ${1\over 4}-{{\epsilon}\over 2}<c<{1\over 4}$,
$$F_2\left(s\right):={1\over {2\pi i}}\int_{\left(-c\right)}{{\Gamma\left
({1\over 4}\pm it+T\right)\Gamma\left(S+T\right)\Gamma\left(-s+T\right
)\Gamma\left(-T\right)\Gamma\left(D-S+s-T\right)}\over {\Gamma\left
({1\over 2}\pm it_1+T\right)}}dT.$$
Using these identities and the definition of $E^{+}\left(S,D\right
)$,
$E^{-}\left(S,D\right)$, $C_{n,\chi}\left(t_1\right)$, (3.6) and that (3.6) is symmetric in $
a$, $b$, $c$, we get for $D$ with large enough real part that (5.29) equals
$${{4\pi\Gamma\left({1\over 2}\pm it_1\right)}\over {2\pi i}}\int_{\left
({1\over 2}-\epsilon\right)}\int_{-\infty}^{\infty}e^{\delta S^2}\left
(4\pi^2\right)^{-S}{{\zeta\left(2S\right)L\left(S\right)}\over {\Gamma\left
(D-S\right)}}{{\Gamma\left({1\over 4}\pm it\pm it_1\right)}\over {
\Gamma\left(\pm 2it\right)}}\chi (t)M\left(S,t\right)dtdS,\eqno (
5.30)$$
where
$$M\left(S,t\right):=M_1\left(S,t\right)+M_2\left(S,t\right),$$
and $M_1\left(S,t\right)$ denotes
$${1\over {2\pi i}}\int_{\left(-{1\over 4}-{{\epsilon}\over 2}\right
)}{{\Gamma\left(-{1\over 2}\pm it_1+S-s\right)\Gamma\left({1\over
2}\pm it_2+s\right)}\over {\Gamma\left(1+s\right)\Gamma\left(D+S-
s-1\right)}}F_1(s)ds,$$
$M_2\left(S,t\right)$ denotes
$${{\sin\pi s_2}\over {\sin\pi s_1}}{1\over {2\pi i}}\int_{\left(
-{1\over 4}-{{\epsilon}\over 2}\right)}{{\Gamma\left(-{1\over 2}\pm
it_2+S-s\right)\Gamma\left({1\over 2}\pm it_1+s\right)\Gamma\left
(1+s-S\right)}\over {\Gamma\left(1+s\right)\Gamma\left(-s\right)\Gamma\left
(D+s\right)}}F_2(s)ds.$$
One can check that (5.30) is a regular function of $D$ for
$\hbox{\rm Re$\, D\ge{1\over 2}$}$, hence by analytic continuation
this equals (5.29) also for $D={1\over 2}$. In the case $D={1\over
2}$ we can
apply Lemma 6.2 to determine $M\left(S,t\right)$. Hence
we proved for any $\delta >0$ that in the case $D={1\over 2}$
(5.29) equals (5.30) with $M\left(S,t\right)$ given by the sum of
(6.21) and (6.22). Recalling the definition of $N\left(S,t\right)$ and $
H_{\chi}\left(S\right)$
from the Introduction we see that (5.29) for $D={1\over 2}$ equals
$$-{{4\pi\Gamma\left({1\over 2}\pm it_1\right)}\over {2\pi i}}{{\sin
\pi s_2}\over {\sin\pi s_1}}\int_{\left({1\over 2}-\epsilon\right
)}e^{\delta S^2}\left(4\pi^2\right)^{-S}\zeta\left(2S\right)L\left
(S\right)\Gamma\left(S\right)\Gamma\left(1-S\right)H_{\chi}\left(
S\right)dS.\eqno (5.31)$$
Assume that $\beta$ in Theorem 1.1 is large enough. Applying Lemma 6.3 (ii) and a convexity bound we see that (5.28) equals (5.31) by writing
$\delta =0$ there. Using Lemma 6.3 (ii) again we see that we can shift
the line of integration to $\hbox{\rm Re$\, S={1\over 2}$ }$in
(5.31). Hence we proved finally that (5.25) equals
$$-4\pi\Gamma\left({1\over 2}\pm it_1\right){{\sin\pi s_2}\over {\sin
\pi s_1}}{1\over {2\pi i}}\int_{\left({1\over 2}\right)}\left(4\pi^
2\right)^{-S}\Gamma\left(S\right)\Gamma\left(1-S\right)\zeta\left
(2S\right)L\left(S\right)H_{\chi}\left(S\right)dS.$$
Using this last relation, (5.1), (5.2), (5.3), (5.20), (5.21),
(5.22), (5.25), (5.26) and (5.27), taking into account the definition of $
C_{0,\chi}$ in
Lemma 4.2 and $\Gamma\left({1\over 2}\right)=\pi^{1/2}$, we get Theorem 1.1 for the case when $
\chi$ satisfies
Condition $D$.

{\bf 5.2. The general case.} To extend the theorem for the general case, we first
need a lemma.

{\bf LEMMA 5.1.} {\it Let} $\beta >0$ {\it and let} $\chi$ {\it be an even holomorphic function on the strip} $\left
|\hbox{\rm Im}\,z\right|<\beta$ {\it such that for a fixed} $A>0$ {\it the function} $\left
|\chi (z)\right|e^{A\left|z\right|^2}$ {\it is bounded}
{\it on the strip} $\left|\hbox{\rm Im}\,z\right|<\beta$. {\it Then for every} $
0<\gamma <\beta$ {\it there is a sequence} $\chi_n$ {\it of entire functions, and a nonnegative function} $
M$ {\it on  }
$[0,\infty )$ {\it with the following properties. The function} $
\chi_n$
{\it satisfies Condition D for every} $n$, {\it for every fixed} $
K>0$
{\it the function} $M(R)e^{KR}$ {\it is bounded on} $[0,\infty )${\it , we have }
$\left|\chi_n\left(z\right)\right|\le M\left(\left|z\right|\right
)$ {\it for every} $n\ge 1$ {\it and} $\left|\hbox{\rm Im}\,z\right
|<\gamma${\it , and finally,} $\chi_n(z)\rightarrow\chi (z)$ {\it for every} $\left
|\hbox{\rm Im}\,z\right|<\gamma${\it .}

{\it Proof.\/} It follows from elementary facts on Fourier transforms that
$$\chi (z)e^{{A\over 2}z^2}=\int_{-\infty}^{\infty}h(x)e^{ixz}dx$$
for $\left|\hbox{\rm Im}\,z\right|<\beta$ where $h$ is an even function such that $
h(x)\ll_{\delta}e^{-\delta\left|x\right|}$
for every $0<\delta <\beta$. Define now
$$\chi_n(z):=e^{-{A\over 2}z^2}\int_{-n}^nh(x)e^{ixz}dx,$$
then for $\left|\hbox{\rm Im}\,z\right|<\gamma$ we have
$$\left|\chi_n\left(z\right)\right|\le\left|e^{-{A\over 2}z^2}\right
|\int_{-\infty}^{\infty}\left|h\left(x\right)\right|e^{\gamma\left
|x\right|}dx.$$
The lemma follows.

Note that using the convexity bound we see that there is a constant $
\beta_0>0$ such that
$${1\over {2\pi i}}\int_{\left({1\over 2}\right)}\left|\zeta\left
(2S\right)L\left(S\right)\right|\left|S\right|^{-{1\over 2}-2\beta_
0}dS<\infty .$$
We choose $\beta$ such that $\beta >\beta_0$. Let $\chi$ be a function satisfying Condition $
C_{\beta}$. Then the
sum in (1.7) and the integral in (1.8) are absolutely
convergent by [B3], formulas (5.2) and (5.3). Then it
follows from Lemma 6.3 (ii) and the dominated
convergence theorem that it is enough to prove Theorem
1.1 for every function $\chi (z)e^{-z^2/N}$ ($N$ is a positive integer) instead of $
\chi$. So we may assume that there is an $A>0$ such that $\chi (z
)e^{A\left|z\right|^2}$
is bounded on the strip $\left|\hbox{\rm Im$\,z$}\right|<\beta$. Finally, for such
functions the theorem follows from Lemma 5.1, Lemma 6.3
(ii), the dominated convergence theorem and the already
proved special case of Theorem 1.1. The theorem is proved.

\noindent{\bf 6. On the kernel function and the integral transform }
\medskip

In this section $t_1$ and $t_2$ are fixed nonzero real numbers, and we
write $s_j={1\over 2}+it$$_j$ for $j=1,2$.

{\bf 6.1. Determination of the kernel function.} The first lemma is proved here in a slightly more general form than necessary; in fact, for Theorem 1.1 we use only the
$n=0$ case. The $n\ge 1$ case would be needed for the proof
of Theorem 1.2. Our main result in this subsection is
Lemma 6.2.

{\bf LEMMA 6.1.} {\it Let} $\epsilon ={1\over {100}}$, {\it and let} $
S$, $B$ {\it and integer} $n$ {\it be given such that} $\hbox{\rm Re$\,
S={1\over 2}-\epsilon$}$, {\it and either}
$$\hbox{\rm $n=0$},\qquad{3\over 4}-{{\epsilon}\over 2}<\hbox{\rm Re$\,
B$}<{3\over 4},$$
{\it or}
$$\hbox{\rm $n\ge 1$},\qquad\hbox{\rm $B$}={1\over 2}.$$
{\it Let} $\gamma_1$ {\it and} $\gamma_2$ {\it be curves (in} $s${\it ) connecting} $
-i\infty$ {\it and} $i\infty$ {\it such that }
$$\hbox{\rm the poles of }\Gamma\left({1\over 2}\pm it_2+s\right)
\Gamma\left(1-n+s-S\right)\hbox{\rm \ lie to the left of $\gamma_
1$},$$
$$\hbox{\rm the poles of }\Gamma\left(-{1\over 2}\pm it_1+S-s\right
)\Gamma\left(n-1+B-s\right)\hbox{\rm \ lie to the right of $\gamma_
1$},$$
$$\hbox{\rm the poles of }\Gamma\left({1\over 2}\pm it_1+s\right)
\Gamma\left(1+n+s-S\right)\Gamma\left(B-S+s+n\right)\hbox{\rm \ lie to the left of $
\gamma_2$},$$
$$\hbox{\rm the poles of }\Gamma\left(-{1\over 2}\pm it_2+S-s\right
)\hbox{\rm \ lie to the right of $\gamma_2$}.$$
{\it In the case} $n=0$ {\it both} $ ${\it of} $\gamma_1$ {\it and} $
\gamma_2$ {\it
may be the line with real part} $-{1\over 4}-{{\epsilon}\over 2}${\it .}

{\it Consider the integrals}
$${1\over {2\pi i}}\int_{\gamma_1}{{\Gamma\left(-{1\over 2}\pm it_
1+S-s\right)\Gamma\left({1\over 2}\pm it_2+s\right)\Gamma\left(1-
n+s-S\right)\Gamma\left(n-1+B-s\right)}\over {\Gamma\left(1-n+s\right
)\Gamma\left(n-1+B-s+S\right)}}ds\eqno (6.1)$$
{\it and}
$${{{{\sin\pi s_2}\over {\sin\pi s_1}}}\over {2\pi i}}\int_{\gamma_
2}{{\Gamma\left(-{1\over 2}\pm it_2+S-s\right)\Gamma\left({1\over
2}\pm it_1+s\right)\Gamma\left(1+n+s-S\right)\Gamma\left(B-S+s+n\right
)}\over {\Gamma\left(B+n+s\right)\Gamma\left(1+n+s\right)}}ds.\eqno
(6.2)$$
{\it Then (6.1) equals}
$$\left(-1\right)^{n-1}\Gamma\left(B-S\right)\Gamma\left(1-S\right
)\Gamma\left({1\over 2}-n\pm it_1\right){{\sin\pi s_2}\over {\sin
\pi s_1}}\left(C_1^{+}Q^{+}+C_1^{-}Q^{-}\right),\eqno (6.3)$$
{\it and (6.2) equals}
$$\left(-1\right)^{n-1}\Gamma\left(B-S\right)\Gamma\left(1-S\right
)\Gamma\left({1\over 2}-n\pm it_1\right){{\sin\pi s_2}\over {\sin
\pi s_1}}\left(C_2^{+}Q^{+}+C_2^{-}Q^{-}\right),\eqno (6.4)$$
{\it where}
$$C_1^{+}:={{\Gamma\left(B-{1\over 2}+n+it_2\right)\Gamma\left({1\over
2}+n+it_2\right)\Gamma\left(S\pm it_1+it_2\right)}\over {\sin\pi\left
(2it_2\right)}}\sin\pi s_1,\eqno $$
$$C_1^{-}:={{\Gamma\left(B-{1\over 2}+n-it_2\right)\Gamma\left({1\over
2}+n-it_2\right)\Gamma\left(S\pm it_1-it_2\right)}\over {\sin\pi\left
(-2it_2\right)}}\sin\pi s_1,$$
$$C_2^{+}:={{\Gamma\left(B-{1\over 2}+n+it_2\right)\Gamma\left({1\over
2}+n+it_2\right)\Gamma\left(S\pm it_1+it_2\right)}\over {\sin\pi\left
(2it_2\right)}}\sin\pi\left({1\over 2}-it_2-S\right),$$
$$C_2^{-}:={{\Gamma\left(B-{1\over 2}+n-it_2\right)\Gamma\left({1\over
2}+n-it_2\right)\Gamma\left(S\pm it_1-it_2\right)}\over {\sin\pi\left
(-2it_2\right)}}\sin\pi\left({1\over 2}+it_2-S\right),$$
$$Q^{+}:=\phi_{i\left({1\over 2}-S\right)}\left(i\left({{1-B}\over
2}-n\right);1-{B\over 2}+it_2,{B\over 2}+it_1,{B\over 2}-it_1,1-{
B\over 2}-it_2\right),\eqno (6.5)$$
$$Q^{-}:=\phi_{i\left({1\over 2}-S\right)}\left(i\left({{1-B}\over
2}-n\right);1-{B\over 2}-it_2,{B\over 2}+it_1,{B\over 2}-it_1,1-{
B\over 2}+it_2\right).\eqno (6.6)$$
{\it Proof.\/} Formula (6.1) equals, by shifting the integration to the left, the sum of
$${{\Gamma\left(1+2it_2\right)\Gamma\left(-2it_2\right)\Gamma\left
({1\over 2}-it_2-n-S\right)\Gamma\left({1\over 2}+n+it_2+S\right)}\over {
\Gamma\left({1\over 2}-it_2-n\right)\Gamma\left({1\over 2}+it_2+n\right
)}}F^{+},\eqno (6.7)$$
$${{\Gamma\left(1-2it_2\right)\Gamma\left(2it_2\right)\Gamma\left
({1\over 2}+it_2-n-S\right)\Gamma\left({1\over 2}+n-it_2+S\right)}\over {
\Gamma\left({1\over 2}+it_2-n\right)\Gamma\left({1\over 2}-it_2+n\right
)}}F^{-}\eqno (6.8)$$
and
$${{\Gamma\left({3\over 2}-n\pm it_2-S\right)\Gamma\left(-{1\over
2}\pm it_2+n+S\right)}\over {\Gamma\left(S\right)\Gamma\left(1-S\right
)}}G,\eqno (6.9)$$
where we write
$$F^{+}:=\sum_{m=0}^{\infty}{{\Gamma\left(B-{1\over 2}+n+it_2+m\right
)\Gamma\left({1\over 2}+n+it_2+m\right)\Gamma\left(\pm it_1+S+it_
2+m\right)}\over {m!\Gamma\left(1+2it_2+m\right)\Gamma\left(B-{1\over
2}+n+it_2+S+m\right)\Gamma\left({1\over 2}+n+it_2+S+m\right)}},$$
$$F^{-}:=\sum_{m=0}^{\infty}{{\Gamma\left(B-{1\over 2}+n-it_2+m\right
)\Gamma\left({1\over 2}+n-it_2+m\right)\Gamma\left(\pm it_1+S-it_
2+m\right)}\over {m!\Gamma\left(1-2it_2+m\right)\Gamma\left(B-{1\over
2}+n-it_2+S+m\right)\Gamma\left({1\over 2}+n-it_2+S+m\right)}},$$
$$G:=\sum_{m=0}^{\infty}{{\Gamma\left({1\over 2}\pm it_1-n+m\right
)\Gamma\left(B-S+m\right)\Gamma\left(1-S+m\right)}\over {m!\Gamma\left
({3\over 2}-n\pm it_2-S+m\right)\Gamma\left(B+m\right)}}.$$
By (3.31), (3.32), (3.33) and (3.34) we have that
$$F^{+}-F^{-}={{\Gamma\left(B-{1\over 2}+n\pm it_2\right)\Gamma\left
({1\over 2}+n\pm it_2\right)\Gamma\left(\pm it_1+S\pm it_2\right)}\over {
\Gamma\left(1+2it_2\right)\Gamma\left(-2it_2\right)}}P,\eqno (6.1
0)$$
and that $F^{+}-G$ equals
$${{\Gamma\left(B-{1\over 2}+n+it_2\right)\Gamma\left({1\over 2}+
n+it_2\right)\Gamma\left(S\pm it_1+it_2\right)\Gamma\left({1\over
2}\pm it_1-n\right)\Gamma\left(B-S\right)\Gamma\left(1-S\right)}\over {
\Gamma\left({1\over 2}+n+it_2+S\right)\Gamma\left({1\over 2}-n-it_
2-S\right)}}\eqno (6.11)$$
times $Q^{+}$, where (the function $\psi$ is defined in (3.28))
$$P:=\psi\left(\alpha ;B-{1\over 2}+n-it_1,{1\over 2}+n-it_1,B-{1\over
2}+n+it_2,{1\over 2}+n+it_2,S+it_2-it_1\right)$$
with the abbreviation
$$\alpha :=B-1+2n+S+it_2-it_1,$$
and
$$Q^{+}:=\psi\left(S+2it_2;S,{3\over 2}+it_2-B-n,S+it_2+it_1,{1\over
2}+n+it_2,S+it_2-it_1\right).\eqno (6.12)$$
Let
$$Q^{-}:=\psi\left(S-2it_2;S,{3\over 2}-it_2-B-n,S-it_2+it_1,{1\over
2}+n-it_2,S-it_2-it_1\right),\eqno (6.13)$$
then by (3.35), (3.36) and (3.37) we have that
$${{\sin\pi\left({1\over 2}+it_2+n+S\right)Q^{+}}\over {\Gamma\left
({1\over 2}+n-it_2\right)\Gamma\left(S-it_2\pm it_1\right)\Gamma\left
(B-{1\over 2}+n-it_2\right)}}+{{\sin\pi\left(2it_2\right)P}\over {
\Gamma\left(B-S\right)\Gamma\left(1-S\right)\Gamma\left({1\over 2}
-n\pm it_1\right)}}\eqno (6.14)$$
 equals
$${{\sin\pi\left({1\over 2}-it_2+n+S\right)Q^{-}}\over {\Gamma\left
({1\over 2}+n+it_2\right)\Gamma\left(S+it_2\pm it_1\right)\Gamma\left
(B-{1\over 2}+n+it_2\right)}}.\eqno (6.15)$$
The identity
$$1={{\sin\pi\left({1\over 2}-it_2-n-S\right)}\over {\sin\pi\left
({1\over 2}+it_2-n-S\right)}}+{{\sin\pi S\sin\pi\left(1+2it_2\right
)}\over {\sin\pi\left({1\over 2}-it_2-n\right)\sin\pi\left({3\over
2}+it_2-n-S\right)}}\eqno (6.16)$$
follows from the easily checked fact that the right-hand
side is a bounded entire function of $S$, and its value is 1 at
$S=0$. Multiplying (6.7) by the right-hand side of (6.16), we see
that the sum of (6.7), (6.8) and (6.9) equals the sum of
$${{\Gamma\left(1-2it_2\right)\Gamma\left(2it_2\right)\Gamma\left
({1\over 2}+it_2-n-S\right)\Gamma\left({1\over 2}+n-it_2+S\right)}\over {
\Gamma\left({1\over 2}+it_2-n\right)\Gamma\left({1\over 2}+n-it_2\right
)}}\left(F^{-}-F^{+}\right)$$
and
$${{\Gamma\left({3\over 2}-n\pm it_2-S\right)\Gamma\left(-{1\over
2}\pm it_2+n+S\right)}\over {\Gamma\left(S\right)\Gamma\left(1-S\right
)}}\left(G-F^{+}\right),$$
which sum, by (6.10) and (6.11), equals the sum of
$${{\pi\Gamma\left(B-{1\over 2}+n\pm it_2\right)\Gamma\left({1\over
2}+n+it_2\right)\Gamma\left(\pm it_1+S\pm it_2\right)}\over {\Gamma\left
({1\over 2}-n+it_2\right)\sin\pi\left({1\over 2}+it_2-n-S\right)}}
P\eqno (6.17)$$
and
$${{\pi\Gamma\left(B-{1\over 2}+n+it_2\right)\Gamma\left({1\over
2}+n+it_2\right)\Gamma\left(S\pm it_1+it_2\right)\Gamma\left({1\over
2}\pm it_1-n\right)\Gamma\left(B-S\right)}\over {\Gamma\left(S\right
)\sin\pi\left({3\over 2}+it_2-n-S\right)}}Q^{+}.\eqno (6.18)$$
Hence we proved that (6.1) equals the sum of (6.17) and
(6.18).

By shifting the integration to the right, we see that (6.2) equals
$${{\sin\pi s_2}\over {\sin\pi s_1}}\left(\Gamma\left(-2it_2\right
)\Gamma\left(1+2it_2\right)F^{+}+\Gamma\left(2it_2\right)\Gamma\left
(1-2it_2\right)F^{-}\right),$$
which, by (6.10), equals
$${{\sin\pi s_2}\over {\sin\pi s_1}}\Gamma\left(B-{1\over 2}+n\pm
it_2\right)\Gamma\left({1\over 2}+n\pm it_2\right)\Gamma\left(\pm
it_1+S\pm it_2\right)P.$$
Hence both (6.1) and (6.2) are linear combinations of $P$ and
$Q^{+}$. By the equality of (6.14) and
(6.15) we can express $P$ by $Q^{+}$ and $Q^{-}$, and by a
tedious, but straightforward calculation we get (6.3) and
(6.4), with $C_1^{+}$, $C_2^{+}$, $C_1^{-}$, $C_2^{-}$ given in the text of the lemma, and $
Q^{+}$,
$Q^{-}$ given by (6.12) and (6.13). During the
calculation we need the identity
$${{\sin\pi S\sin2\pi it_2}\over {\sin\pi\left({1\over 2}+it_2-S\right
)\sin\pi\left({1\over 2}+it_2+S\right)}}+{{\sin\pi s_2}\over {\sin
\pi\left({1\over 2}+it_2-S\right)}}={{\sin\pi s_2}\over {\sin\pi\left
({1\over 2}-it_2-S\right)}};$$
for its proof it is enough to show that the difference
of the two sides is a regular function of $S$, and it is
not hard to see.

By (3.30), (3.23) and (3.24) we get the expressions (6.5) and (6.6) for $
Q^{+}$ and $Q^{-}$. The lemma is proved.

{\bf LEMMA 6.2.} {\it Let} $\epsilon ={1\over {100}}${\it , and let} $
t$ {\it and} $S$ {\it be given such that} $t$ {\it is real and} $\hbox{\rm Re$\,
S={1\over 2}-\epsilon$}$. {\it Consider the integrals}
$${1\over {2\pi i}}\int_{\left(-{1\over 4}-{{\epsilon}\over 2}\right
)}{{\Gamma\left(-{1\over 2}\pm it_1+S-s\right)\Gamma\left({1\over
2}\pm it_2+s\right)}\over {\Gamma\left(1+s\right)\Gamma\left(S-s-{
1\over 2}\right)}}F_1(s)ds\eqno (6.19)$$
 {\it and}
$${{\sin\pi s_2}\over {\sin\pi s_1}}{1\over {2\pi i}}\int_{\left(
-{1\over 4}-{{\epsilon}\over 2}\right)}{{\Gamma\left(-{1\over 2}\pm
it_2+S-s\right)\Gamma\left({1\over 2}\pm it_1+s\right)\Gamma\left
(1+s-S\right)}\over {\Gamma\left(1+s\right)\Gamma\left(-s\right)\Gamma\left
({1\over 2}+s\right)}}F_2(s)ds,\eqno (6.20)$$
{\it where} $F_1(s)$ {\it
denotes}
$${1\over {2\pi i}}\int_{\left(-c\right)}{{\Gamma\left({1\over 4}
\pm it+T\right)\Gamma\left(S+T\right)\Gamma\left(-T\right)\Gamma\left
(1-S+s+T\right)\Gamma\left(-{1\over 2}-s-T\right)}\over {\Gamma\left
({1\over 2}\pm it_1+T\right)}}dT,$$
{\it and} $F_2(s)$ {\it denotes}
$${1\over {2\pi i}}\int_{\left(-c\right)}{{\Gamma\left({1\over 4}
\pm it+T\right)\Gamma\left(S+T\right)\Gamma\left(-T\right)\Gamma\left
(-s+T\right)\Gamma\left({1\over 2}-S+s-T\right)}\over {\Gamma\left
({1\over 2}\pm it_1+T\right)}}dT$$
{\it with} ${1\over 4}-{{\epsilon}\over 2}<c<{1\over 4}$.

{\it Then (6.19) equals}
$$-\Gamma\left({1\over 4}\pm it\right)\Gamma\left({1\over 2}-S\right
)\Gamma\left(S\right)\Gamma\left(1-S\right)\left(A^{+}\left(S,t\right
)+A^{-}\left(S,t\right)\right)\sin\pi s_2\eqno (6.21)$$
{\it with}
$$A^{+}\left(S,t\right):={{\Gamma\left(S\pm it_1+it_2\right)\Gamma\left
({1\over 4}+it_2\pm it\right)}\over {\sin\pi\left(2it_2\right)}}\phi_{
i\left({1\over 2}-S\right)}^{+}\left(t\right),$$

$$A^{-}\left(S,t\right):={{\Gamma\left(S\pm it_1-it_2\right)\Gamma\left
({1\over 4}-it_2\pm it\right)}\over {\sin\pi\left(-2it_2\right)}}
\phi_{i\left({1\over 2}-S\right)}^{-}\left(t\right),$$
{\it and (6.20) equals}
$$-\Gamma\left({1\over 4}\pm it\right)\Gamma\left({1\over 2}-S\right
)\Gamma\left(S\right)\Gamma\left(1-S\right){{\sin\pi s_2}\over {\sin
\pi s_1}}\left(B^{+}\left(S,t\right)+B^{-}\left(S,t\right)\right)\eqno
(6.22)$$
{\it with}
$$B^{+}\left(S,t\right):={{\Gamma\left(S\pm it_1+it_2\right)\Gamma\left
({1\over 4}+it_2\pm it\right)}\over {\sin\pi\left(2it_2\right)}}\left
(\sin\pi\left({1\over 2}-it_2-S\right)\right)\phi_{i\left({1\over
2}-S\right)}^{+}\left(t\right),$$
$$B^{-}\left(S,t\right):={{\Gamma\left(S\pm it_1-it_2\right)\Gamma\left
({1\over 4}-it_2\pm it\right)}\over {\sin\pi\left(-2it_2\right)}}\left
(\sin\pi\left({1\over 2}+it_2-S\right)\right)\phi_{i\left({1\over
2}-S\right)}^{-}\left(t\right).$$
{\it Proof.\/} We see by (3.26) and (3.27) that $F_1(s)$ equals
$$\Gamma\left({1\over 4}\pm it\right)\Gamma\left(-{1\over 4}-s\pm
it\right)\Gamma\left(S\right)\Gamma\left(1+s-S\right)\Gamma\left({
1\over 2}-S\right)\Gamma\left(S-s-{1\over 2}\right)$$
times
$$\phi_{i\left({{1+s}\over 2}-S\right)}\left(t;{1\over 4},{1\over
4}+it_1,{1\over 4}-it_1,{5\over 4}+s\right).\eqno (6.23)$$
Similarly, we see that
$${1\over {2\pi i}}\int_{\left(-{{\epsilon}\over 4}\right)}{{\Gamma\left
({1\over 4}\pm it_1+A\right)\Gamma\left({1\over 4}+A\right)\Gamma\left
(-{1\over 4}+A-s\right)\Gamma\left(-A\pm it\right)}\over {\Gamma\left
({3\over 4}+A-S\right)\Gamma\left(-{1\over 4}+A+S-s\right)}}dA\eqno
(6.24)$$
 equals
$$\Gamma\left({1\over 4}\pm it\pm it_1\right)\Gamma\left({1\over
4}\pm it\right)\Gamma\left(-{1\over 4}\pm it-s\right)$$
times
$$\phi_{i\left({1\over 4}+{s\over 2}\right)}\left(t_1;{1\over 4}+
it,{1\over 2}-S,-{1\over 2}+S-s,{3\over 4}+it\right).\eqno (6.25)$$
We see by (3.25) and by the symmetry of the Wilson
function in its parameters (see the sentence above (3.25))
that (6.23) equals (6.25). Hence $F_1(s)$ equals
$${{\Gamma\left(S\right)\Gamma\left(1+s-S\right)\Gamma\left({1\over
2}-S\right)\Gamma\left(S-s-{1\over 2}\right)}\over {\Gamma\left({
1\over 4}\pm it\pm it_1\right)}}$$
times (6.24). This means that (6.19) equals
$${1\over {2\pi i}}{{\Gamma\left(S\right)\Gamma\left({1\over 2}-S\right
)}\over {\Gamma\left({1\over 4}\pm it\pm it_1\right)}}\int_{\left
(-{{\epsilon}\over 4}\right)}{{\Gamma\left({1\over 4}\pm it_1+A\right
)\Gamma\left({1\over 4}+A\right)\Gamma\left(-A\pm it\right)}\over {
\Gamma\left({3\over 4}+A-S\right)}}L_1(A)dA,\eqno (6.26)$$
where $L_1(A)$ denotes
$${1\over {2\pi i}}\int_{\left(-{1\over 4}-{{\epsilon}\over 2}\right
)}{{\Gamma\left(-{1\over 2}\pm it_1+S-s\right)\Gamma\left({1\over
2}\pm it_2+s\right)\Gamma\left(1+s-S\right)\Gamma\left(-{1\over 4}
+A-s\right)}\over {\Gamma\left(1+s\right)\Gamma\left(-{1\over 4}+
A-s+S\right)}}ds.\eqno $$
We see by (3.26) and (3.27) that $F_2(s)$ equals
$$\Gamma\left({1\over 4}\pm it\right)\Gamma\left({3\over 4}-S+s\pm
it\right)\Gamma\left(S\right)\Gamma\left(-s\right)\Gamma\left({1\over
2}+s\right)\Gamma\left({1\over 2}-S\right)$$
times
$$\phi_{i\left({{s+S}\over 2}\right)}\left(t;{1\over 4},{1\over 4}
+it_1,{1\over 4}-it_1,{1\over 4}+S-s\right).\eqno (6.27)$$
Similarly, we see that
$${1\over {2\pi i}}\int_{\left(-{{\epsilon}\over 4}\right)}{{\Gamma\left
({1\over 4}\pm it_1+A\right)\Gamma\left({1\over 4}+A\right)\Gamma\left
({3\over 4}+A-S+s\right)\Gamma\left(-A\pm it\right)}\over {\Gamma\left
({3\over 4}+A-S\right)\Gamma\left({3\over 4}+A+s\right)}}dA\eqno
(6.28)$$
equals
$$\Gamma\left({1\over 4}\pm it\pm it_1\right)\Gamma\left({1\over
4}\pm it\right)\Gamma\left({3\over 4}\pm it+s-S\right)$$
times
$$\phi_{i\left({1\over 4}+{{s-S}\over 2}\right)}\left(t_1;{1\over
4}+it,{1\over 2}-S,{1\over 2}+s,{3\over 4}+it\right).\eqno (6.29)$$
We see again by (3.25) and by the symmetry of the Wilson
function in its parameters that (6.27) equals (6.29). Hence $F_2(
s)$ equals
$${{\Gamma\left(S\right)\Gamma\left(-s\right)\Gamma\left({1\over
2}+s\right)\Gamma\left({1\over 2}-S\right)}\over {\Gamma\left({1\over
4}\pm it\pm it_1\right)}}$$
times (6.28). This means that (6.20) equals
$${1\over {2\pi i}}{{\sin\pi s_2}\over {\sin\pi s_1}}{{\Gamma\left
(S\right)\Gamma\left({1\over 2}-S\right)}\over {\Gamma\left({1\over
4}\pm it\pm it_1\right)}}\int_{\left(-{{\epsilon}\over 4}\right)}{{
\Gamma\left({1\over 4}\pm it_1+A\right)\Gamma\left({1\over 4}+A\right
)\Gamma\left(-A\pm it\right)}\over {\Gamma\left({3\over 4}+A-S\right
)}}L_2(A)dA,\eqno (6.30)$$
where $L_2(A)$ denotes
$${1\over {2\pi i}}\int_{\left(-{1\over 4}-{{\epsilon}\over 2}\right
)}{{\Gamma\left(-{1\over 2}\pm it_2+S-s\right)\Gamma\left({1\over
2}\pm it_1+s\right)\Gamma\left(1+s-S\right)\Gamma\left({3\over 4}
+A-S+s\right)}\over {\Gamma\left({3\over 4}+A+s\right)\Gamma\left
(1+s\right)}}ds.$$
Applying Lemma 6.1 with $n=0,$ $B=A+{3\over 4}$ we see for $\hbox{\rm Re$\,
A=-{{\epsilon}\over 4}$}$
that $L_1(A)$ equals (6.3), and ${{\sin\pi s_2}\over {\sin\pi s_1}}
L_2(A)$ equals (6.4) with $n=0,$ $B=A+{3\over 4}$ there.

We see by (6.5), (3.26) and (3.27) (using again that the Wilson function is
symmetric in the parameters $a,b,c$ and $1-d$) that if $n=0,$ $B=
A+{3\over 4}$, $\hbox{\rm Re$\, A=-{{\epsilon}\over 4}$}$, then $Q^{
+}$ equals
$${1\over {\Gamma\left({1\over 2}\pm ix\right)\Gamma\left({1\over
2}\pm ix-it_1+it_2\right)\Gamma\left({1\over 2}+it_1\right)\Gamma\left
({1\over 2}+it_2\right)\Gamma\left({1\over 4}+A+it_1\right)\Gamma\left
({1\over 4}+A+it_2\right)}}$$
times
$${1\over {2\pi i}}\int_{\left(d\right)}{{\Gamma\left({1\over 2}\pm
ix+R\right)\Gamma\left({1\over 2}+it_1+R\right)\Gamma\left({1\over
4}+A+it_1+R\right)\Gamma\left(-R\right)\Gamma\left(it_2-it_1-R\right
)}\over {\Gamma\left({3\over 4}+A+R\right)\Gamma\left(1+it_1+it_2
+R\right)}}dR$$
with $d=-1/8$, where we write $S={1\over 2}+ix$.

Observe that
$${1\over {2\pi i}}\int_{\left(-{{\epsilon}\over 4}\right)}{{\Gamma\left
({1\over 4}-it_1+A\right)\Gamma\left({1\over 4}+A\right)\Gamma\left
(-A\pm it\right)\Gamma\left({1\over 4}+A+it_1+R\right)}\over {\Gamma\left
({3\over 4}+A+R\right)}}dA\eqno (6.31)$$
equals
$${{\Gamma\left({1\over 4}-it_1\pm it\right)\Gamma\left({1\over 4}
\pm it\right)\Gamma\left({1\over 4}+it_1+R\pm it\right)}\over {\Gamma\left
({1\over 2}+R\right)\Gamma\left({1\over 2}+it_1+R\right)\Gamma\left
({1\over 2}-it_1\right)}}\eqno (6.32)$$
by (3.16) and (3.17).

We claim that (6.19) equals
$$-{{\Gamma\left(S+it_1+it_2\right)\Gamma\left({1\over 4}\pm it\right
)\Gamma\left({1\over 2}-S\right)}\over {\Gamma\left(1-S-it_1+it_2\right
)\sin\pi\left(2it_2\right)\Gamma\left({1\over 4}\pm it+it_1\right
)}}\sin\pi s_2\eqno (6.33)$$
times
$${1\over {2\pi i}}\int_{\left(-{1\over 8}\right)}{{\Gamma\left(S
+R\right)\Gamma\left(1-S+R\right)\Gamma\left({1\over 4}+it_1+R\pm
it\right)\Gamma\left(-R\right)\Gamma\left(-it_1+it_2-R\right)}\over {
\Gamma\left(1+it_1+it_2+R\right)\Gamma\left({1\over 2}+R\right)}}
dR\eqno (6.34)$$
plus the similar product obtained by writing $-t_2$ in place of $
t_2$ in
(6.33) and (6.34). Indeed, we can see it by (6.26), (6.31), (6.32), by the above-mentioned fact that $
L_1(A)$
equals (6.3) writing $n=0,$ $B=A+{3\over 4}$ there, by the above expression for $
Q^{+}$, and
by the fact that $Q^{-}$ and $C_1^{-}$ are obtained from $Q^{+}$ and $
C_1^{+}$
by writing $-t_2$ in place of $t_2$.

Similarly, but using (6.30) in place of (6.26), we see that
(6.20) equals
$$-{{\Gamma\left(S+it_1+it_2\right)\Gamma\left({1\over 4}\pm it\right
)\Gamma\left({1\over 2}-S\right)}\over {\Gamma\left(1-S-it_1+it_2\right
)\sin\pi\left(2it_2\right)\Gamma\left({1\over 4}\pm it+it_1\right
)}}{{\sin\pi s_2}\over {\sin\pi s_1}}\sin\pi\left({1\over 2}-it_2
-S\right)\eqno (6.35)$$
times (6.34) plus the similar product obtained by writing
$-t_2$ in place of $t_2$ in (6.35) and (6.34).

By (3.26) and (3.27) we see that (6.34) equals
$$\Gamma\left(S\right)\Gamma\left(1-S\right)\Gamma\left(S-it_1+it_
2\right)\Gamma\left(1-S-it_1+it_2\right)\Gamma\left({1\over 4}+it_
1\pm it\right)\Gamma\left({1\over 4}+it_2\pm it\right)$$

times
$$\phi_{i\left({1\over 2}-S\right)}\left(t;{3\over 4}+it_2,{1\over
4}+it_1,{1\over 4}-it_1,{3\over 4}-it_2\right).$$
This proves the lemma.

{\bf 6.2. An estimation for} $H_{\chi}\left(S\right)${\bf .} During the proof of
Theorem 1.1 we need an upper bound for $H_{\chi}\left(S\right)$ defined in Theorem 1.1 not only for an individual
$\chi$ but also for a function series $\chi_n$, assuming a universal
upper bound for every $\left|\chi_n\right|$. The most important aspect of the lemma below is that the estimate (6.36) depends only on the upper bound $
M$ for $\left|\chi\right|$.

{\bf LEMMA 6.3.} {\it (i) Recall the notations $\phi_{\lambda}^{+}\left(x\right)$, $\phi_{\lambda}^{-}\left(x\right)$ from Section 1.2. There is an absolute constant} $C>0$ {\it such that we have}
$$\left|\phi_{i\left({1\over 2}-S\right)}^{+}\left(t\right)\right
|+\left|\phi_{i\left({1\over 2}-S\right)}^{-}\left(t\right)\right
|\ll e^{\pi\left(\left|S\right|+\left|t\right|\right)}\left(1+\left
|S\right|\right)^C\left(1+\left|t\right|\right)^C$$
{\it with an implied absolute constant for every} $S$ {\it with} $
-1\le\hbox{\rm $ $Re$\, S\le 2$}$ {\it and for every real} $t$.

{\it (ii) Let} $\beta >0$ {\it be a given number and let} $M$ {\it be a given nonnegative function on} $[0,\infty )$ {\it satisfying that for every fixed} $K>0$ {\it the function} $
M(R)e^{-\pi R}\left(1+R\right)^K$ {\it is bounded on} $[0,\infty
)${\it . Then, if }
$\chi$ {\it is any even holomorphic function on the strip} $\left
|\hbox{\rm Im}\,z\right|<\beta$ $wi${\it th} $\left|\chi\left(z\right
)\right|\le M\left(\left|z\right|\right)$ {\it on this strip, then for every} $0<B<{1\over 2}+2\beta$ {\it we have that} $H_{\chi}\left
(S\right)$ {\it is regular in the strip} ${1\over 2}-{1\over {100}}
\le\hbox{\rm Re$\, S\le{1\over 2}$}${\it , and for every} $S$ {\it in this strip we have}
$$\Gamma^2\left(1-S\right)H_{\chi}\left(S\right)\ll_{B,M}\left(1+\left
|S\right|\right)^{-B}.\eqno (6.36)$$
{\it Proof.\/} To show (i) note that for every fixed real $t$ the functions $
\phi_{i\left({1\over 2}-S\right)}^{+}\left(t\right)$ and $\phi_{i\left
({1\over 2}-S\right)}^{-}\left(t\right)$ are
entire in $S$. Combining this fact with (3.26) and (3.27) we
get (i) by trivial estimates. The regularity statement in (ii) follows then at once from (i) and from the definition.

By the definition of $N^{+}\left(S,t\right)$ in Section 1, and by (3.26), (3.27) we see for
any real $t$ and for any $S$ with ${1\over 2}-{1\over {100}}
\le\hbox{\rm Re$\, S\le{1\over 2}$}$ that
$${{\Gamma\left(S\right)\Gamma\left(1-S\right)N^{+}\left(S,t\right
)}\over {\Gamma\left(S+it_1\pm it_2\right)}}$$
equals
$${{\sin\pi\left(S+it_1-it_2\right)}\over {\pi\sin\pi\left(2it_2\right
)\Gamma\left({1\over 4}+it_1\pm it\right)}}\left(\sin\pi s_1+\sin
\pi\left({1\over 2}-it_2-S\right)\right)\eqno (6.37)$$
times
$${1\over {2\pi i}}\int_{\left(-1/8\right)}{{\Gamma\left(S+R\right
)\Gamma\left(1-S+R\right)\Gamma\left({1\over 4}+it_1+R\pm it\right
)\Gamma\left(-R\right)\Gamma\left(-it_1+it_2-R\right)}\over {\Gamma\left
(1+it_1+it_2+R\right)\Gamma\left({1\over 2}+R\right)}}dR.\eqno (6
.38)$$
We get
$$\hbox{\rm }{{\Gamma\left(S\right)\Gamma\left(1-S\right)N^{-}\left
(S,t\right)}\over {\Gamma\left(S+it_1\pm it_2\right)}}$$
by writing $-t_2$ in place of $t_2$ in (6.37) and (6.38).

We now show (6.36). It is clear that taking the term
$\sin\pi s_1$ from the bracket in (6.37) we get expressions
acceptable in (6.36). On the other hand, we have that
$$\sin\pi\left(S+it_1-it_2\right)\sin\pi\left({1\over 2}-it_2-S\right
)={{\cos\pi\left(2S+it_1-{1\over 2}\right)}\over 2}-{{\cos\pi\left
(-2it_2+it_1+{1\over 2}\right)}\over 2}$$
by [G-R], p. 29, 1.313.5. Taking the second term from here
in (6.37) gives again an acceptable contribution in (6.36); the first term is independent of
$t_2$. So defining
$$G(R):=\int_{-\infty}^{\infty}{{\Gamma\left({1\over 4}\pm it\right
)\Gamma\left({1\over 4}-it_1\pm it\right)\Gamma\left({1\over 4}+i
t_1+R\pm it\right)}\over {\Gamma\left(\pm 2it\right)}}\chi (t)dt,\eqno
(6.39)$$
it is enough to prove for ${1\over 2}-{1\over {100}}
\le\hbox{\rm Re$\, S\le{1\over 2}$}$ that the difference of
$${1\over {2\pi i}}\int_{\left(-1/8\right)}{{\Gamma\left(S+R\right
)\Gamma\left(1-S+R\right)G\left(R\right)\Gamma\left(-R\right)\Gamma\left
(-it_1+it_2-R\right)}\over {\Gamma\left(1+it_1+it_2+R\right)\Gamma\left
({1\over 2}+R\right)}}dR$$
and the same integral with $-t_2$ in place of $t_2$ is
$\ll_{B,M}e^{-\pi\left|S\right|}\left(1+\left|S\right|\right)^{-B}
.$ We claim that
$${{\Gamma\left(-it_1+it_2-R\right)\Gamma\left(1+it_1-it_2+R\right
)}\over {\Gamma\left(1+it_1+it_2+R\right)\Gamma\left(-it_1-it_2-R\right
)}}-1$$
equals
$${{\Gamma\left(-it_1+it_2-R\right)\Gamma\left(1+it_1-it_2+R\right
)\Gamma\left({1\over 2}+it_2\right)\Gamma\left({1\over 2}-it_2\right
)}\over {\Gamma\left({1\over 2}+it_1+R\right)\Gamma\left({1\over
2}-it_1-R\right)\Gamma\left(1+2it_2\right)\Gamma\left(-2it_2\right
)}}.$$
This is true because the difference of these two functions is a bounded entire function of $
R$ which vanishes at $R={1\over 2}-it_1$. Using this identity we see that it is enough to prove for
${1\over 2}-{1\over {100}}
\le\hbox{\rm Re$\, S\le{1\over 2}$}$ that
$${1\over {2\pi i}}\int_{\left(-1/8\right)}{{\Gamma\left(S+R\right
)\Gamma\left(1-S+R\right)G\left(R\right)\Gamma\left(-R\right)\Gamma\left
(-it_1\pm it_2-R\right)}\over {\Gamma\left({1\over 2}+it_1+R\right
)\Gamma\left({1\over 2}-it_1-R\right)\Gamma\left({1\over 2}+R\right
)}}dR$$
is $\ll_{B,M}e^{-\pi\left|S\right|}\left(1+\left|S\right|\right)^{-
B}.$
By shifting the $R$-integration to the left, we see then
that it is enough to prove that
$$H\left(R\right):={{G\left(R\right)}\over {\Gamma\left({1\over 2}
+it_1+R\right)\Gamma\left({1\over 2}+R\right)}}$$
is holomorphic for $\hbox{\rm Re$\, R>-{1\over 4}$}-\beta$ and satisfies
$${{G\left(R\right)}\over {\Gamma\left({1\over 2}+it_1+R\right)\Gamma\left
({1\over 2}+R\right)}}\ll_{K,\rho ,M}e^{\pi\left|R\right|}\left(1
+\left|R\right|\right)^{-K}$$
for every $K>0$ and $0\le\rho <\beta$ on the strip
$-{1\over 4}-\rho\le\hbox{\rm Re$\, R$}\le 1.$
We now prove this statement. It is clear that we may
assume that ${1\over 4}+\rho$ and ${1\over 4}-\rho$ are not integers.

Let $b$ be a large positive integer. There are constants $c_{a,b}$
such that
$${{\Gamma\left(z\right)}\over {\Gamma\left(z+b\right)}}=\sum_{a=
0}^{b-1}c_{a,b}\left(z+a\right)^{-1}.$$
Applying it for $z={1\over 4}+it_1+R\pm it$, we see that $\Gamma\left
({1\over 4}+it_1+R\pm it\right)$ equals
$$\Gamma\left({1\over 4}+it_1+R\pm it+b\right)\sum_{0\le a_1,a_2\le
b-1}{{c_{a_1,b}c_{a_2,b}}\over {\left({1\over 4}+it_1+R+it+a_1\right
)\left({1\over 4}+it_1+R-it+a_2\right)}}.$$
We use that
$${{-2it+a_2-a_1}\over {\left({1\over 4}+it_1+R+it+a_1\right)\left
({1\over 4}+it_1+R-it+a_2\right)}}$$
equals
$${1\over {{1\over 4}+it_1+R+it+a_1}}-{1\over {{1\over 4}+it_1+R-
it+a_2}},\eqno (6.40)$$
and because of the presence of the factor
${1\over {\Gamma\left(\pm 2it\right)}}$, shifting the line of integration in (6.39) to $\hbox{\rm Im$
\,t=\pm\rho$}$ (the
minus sign is used in the case of the first term in (6.40), and the
plus sign in the case of the second term), we get such an
expression for $G(R)$ which proves the above statement for
the function $H\left(R\right)$. We cross some poles when we shift the
$t$-integration, but the residues also give holomorphic
expressions for $H\left(R\right)$ in the required strip, because of
the factor $\Gamma\left({1\over 2}+it_1+R\right)\Gamma\left({1\over
2}+R\right)$ in the denominator of
$H\left(R\right)$. The lemma is proved.

\vskip10pt

\centerline{\bf Important notations}

\vskip20pt

$\matrix{\left(a\right)_n&\quad\xref{poch}\cr
\arg z&\quad\xref{arg}\cr
B_0(z)&\quad\xref{B0}\cr
B_n(z)&\quad\xref{Bn}\cr
{\rm C}{\rm o}{\rm n}{\rm d}{\rm i}{\rm t}{\rm i}{\rm o}{\rm n}C_{
\beta}&\quad\xref{Cbeta}\cr
{\rm C}{\rm o}{\rm n}{\rm d}{\rm i}{\rm t}{\rm i}{\rm o}{\rm n}D&
\quad\xref{Cond(D)}\cr
D_1&\quad\xref{D1},\xref{D1konk}\cr
D_4&\quad\xref{D4},\xref{D4konk}\cr
d\mu_z&\quad\xref{dmuz}\cr
\Delta_l&\quad\xref{Laplace}\cr
\delta_{u_1,u_2}&\quad\xref{kron}\cr
\phantom{\Gamma\left(X\pm Y\right),\Gamma\left(X\pm Y\pm Z\right)}&\quad\cr}
$

$\matrix{E_{\eufm{a}}\left(z,s,{1\over 2}\right)&\quad\xref{Eisenstein}\cr
e(x)&\quad\xref{e(x)}\cr
\left(f_1,f_2\right)_1&\quad\xref{f1f2}\cr
\left(f_1,f_2\right)_4&\quad\xref{f1f24}\cr
_{q+1}F_q&\quad\xref{genhyp}\cr
F\left(\alpha ,\beta ,\gamma ;z\right)&\quad\xref{hyp}\cr
\phi_{\lambda}\left(x;a,b,c,d\right)&\quad\xref{wilson}\cr
\phi_{\lambda}^{+}\left(x\right),\phi_{\lambda}^{-}\left(x\right)&\quad\xref{phi+}\cr
g_{k,j}&\quad\xref{gkj}\cr
\gamma_j&\quad\xref{gammaj}\cr
\Gamma\left(X\pm Y\right),\Gamma\left(X\pm Y\pm Z\right)&\quad\xref{gamma}\cr
\Gamma_0(4)&\quad\xref{Gammazero4}\cr
\Gamma_{\infty}&\quad\xref{GammaInfty}\cr
\bbb H&\quad\xref{H}\cr
H(z,w)&\quad\xref{H(z,w)}\cr
H_{\chi}\left(S\right)&\quad\xref{Hchi(S)}\cr
j_{\gamma}(z)&\quad\xref{jgamma}\cr
K_k&\quad\xref{MaassK}\cr
\left(\kappa (u)\right)(z)&\quad\xref{kappa}\cr
\left(\kappa_n(u)\right)(z)&\quad\xref{kappa(n)}\cr
L\left(S\right)=L\left(S,u_1\otimes\overline {u_2}\right)&\quad\xref{Rankin}\cr
L_l^2(D_4)&\quad\xref{Ll(d4)}\cr
L_k&\quad\xref{MaassL}\cr
N\left(S,t\right)&\quad\xref{N(s,t)}\cr
\nu (\gamma )&\quad\xref{nu}\cr
\psi\left(A;B,C,D,E,F\right)&\quad\xref{psi}\cr
R_l(D_4)&\quad\xref{Rl(D4)}\cr
\rho_{f,\eufm{a}}(m)&\quad\xref{rho(f,a,m)}\cr
\rho_{u_1}(m),\rho_{u_2}(m)&\quad\xref{rhou1rhou2}\cr
{\rm S}{\rm h}{\rm i}{\rm m}F&\quad\xref{Shim}\cr
S_{2k+{1\over 2}}&\quad\xref{hol}\cr
\sigma_{\eufm{a}}&\quad\xref{scaling}\cr
S_n\left(x^2;a,b,c\right)&\quad\xref{Hahn}\cr
t_1,t_2&\quad\xref{tj}\cr
T_j&\quad\xref{Tj}\cr
T_z&\quad\xref{Tz}\cr
\theta\left(z\right)&\quad\xref{theta}\cr
u_1,u_2&\quad\xref{u1u2}\cr
u_{j,1/2}&\quad\xref{u(j,1/2)}\cr
w_{a,b,c}\left(x\right)&\quad\xref{weight}\cr
W\left(A;B,C,D,E,F\right)&\quad\xref{w}\cr
W_{\alpha ,\beta}\left(y\right)&\quad\xref{Whittaker}\cr
\zeta\left(S\right)&\quad\xref{zeta}\cr
\zeta_{\eufm{a}}(f,r)&\quad\xref{zeta(a,f,r)}\cr}
$

 \bigskip\noindent {\bf References}

\nobreak
\parindent=12pt
\nobreak

\item{[A-A-R]} G.E. Andrews, R. Askey, R. Roy, {Special
Functions}, {\it Cambridge Univ. Press,\/} 1999

\item{[B1]} A. Bir\'o,  {\it A relation between triple products of weight 0 and weight 1/2 cusp forms}, Israel J. of Math., 182
(2011), 61-101.

\item{[B2]} A. Bir\'o,  {\it An expansion theorem concerning Wilson functions and polynomials}, Acta Math. Hung., 135 (2012), 350-382.

\item{[B3]} A. Bir\'o,  {\it A duality relation for certain triple products of automorphic forms}, Israel J. of Math., 192
(2012), 587-636.

\item{[B4]} A. Bir\'o, {\it Cycle integrals of Maass forms of weight 0 and Fourier coefficients of Maass forms of weight 1/2}, Acta Arithmetica, 94 (2000), 103-152.

\item{[B-E]} P. Borwein, T. Erdelyi, {Polynomials and polynomial inequalities,} {\it Springer}, 2012

\item{[B-K]} J. Buttcane, R. Khan,  {\it On the fourth moment of Hecke-Maass forms and the random wave conjecture}, Compositio Mathematica, 153 (7) (2017), 1479-1511.

\item{[B-M]} E.M. Baruch, Z. Mao,
{\it A generalized Kohnen-Zagier formula for Maass forms}, J. London Math. Soc. (2), 82 (2010), no. 1,
1-16.

\item{[D-I-T]} W. Duke, O. Imamoglu, \'A. T\'oth, {\it  }
{\it Geometric invariants for real
quadratic fields}, Ann. of Math., 184 (3), (2016),
949-990.

\item{[F]} J.D. Fay, {\it Fourier coefficients of the resolvent for a Fuchsian group}, J. Reine Angew. Math., 294 (1977),
143-203.

\item{[G1]} W. Groenevelt, {\it The Wilson function transform}, Int. Math. Res. Not. 2003 (52), (2003), 2779--2817.

\item{[G2]} W. Groenevelt, {\it Wilson function transforms related to Racah coefficients}, Acta Appl. Math. 91 (2), (2006), 133--191.

\item{[G-R]} I.S. Gradshteyn, I.M. Ryzhik, {Table of integrals, series and
products, 6th edition,} {\it Academic Press}, 2000

\item{[H]} D.A. Hejhal, {The Selberg Trace Formula for $PSL(2,{\bf R}
)$,
vol. 2}, Springer, 1983

\item{[H-K]} P. Humphries, R. Khan, {\it $L^p$-bounds for automorphic forms via spectral reciprocity},  arXiv e-prints (2022), arXiv:2208.05613

\item{[I]} H. Iwaniec, {Introduction to the spectral theory of automorphic forms,} {\it Rev. Mat. Iberoamericana},
1995

\item{[I-K]} H. Iwaniec, E. Kowalski, {Analytic Number Theory,} {\it AMS Colloqium Publications, Vol. 53, Providence RI, American Mathematical Society},
2004

\item{[K]} T.H. Koornwinder, {\it A
new proof of a Paley-Wiener type theorem for the Jacobi
transform,} Ark. Mat., 13 (1975), 145-159.

\item{[K-S]} S. Katok, P. Sarnak, {\it Heegner points, cycles and Maass forms,\/} Israel J. of Math., 84 (1993), 193-227.

\item{[Kw]} C-H. Kwan, {\it Spectral moment formulae for} $GL(3)\times
GL(2)$ $L${\it -functions}, arXiv e-prints (2021), arXiv:2112.08568

\item{[L]} E. Lindenstrauss, {\it Invariant measures and arithmetic quantum unique ergodicity,\/} Ann. of Math., 163 (1), (2006), 165-219.

\item{[Nel1]} P. D. Nelson, {\it Subconvex Equidistribution of Cusp Forms: Reduction to Eisenstein Observables}, Duke Mathematical Journal 168 (9) (2019), 1665-1722.

\item{[Nel2]} P. D. Nelson, {\it The spectral decomposition of $|\theta|^2$}, Math. Z. 298 (3-4), (2021), 1425-1447.

\item{[P]} N.V. Proskurin, {\it On general Klosterman
sums (in Russian)}, Zap. Naucn. Sem. LOMI 302 (2003),
107-134.

\item{[R]} W. Roelcke, {\it Das Eigenwertproblem der automorphen Formen in der hyperbolischen Ebene, I.}, Mathematische Annalen 167 (4) (1966), 292-337.

\item{[Sa]} P. Sarnak, {\it Additive number theory and Maass }
{\it forms,\/} in: Chudnovsky, D.V., Chudnovsky, G.V., Cohn, H., Nathatnson, M.B., (eds) Number Theory. Proceedings, New York 1982 (Lecture Notes Math. vol. 1052., pp. 286-309.), {\it Springer}, 1982

\item{[S]} L.J. Slater, {Generalized hypergeometric functions,} {\it Cambridge Univ. Press},
1966

\item{[So]} K. Soundararajan, {\it Quantum unique
ergodicity for $SL(2,{\bf Z})\setminus \bbb H$}, Ann. of Math., 172 (2), (2010), 1529-1538.

\item{[W]} F.J.W. Whipple, {\it Relations between
well-poised hypergeometric series of the type $_7F$$_6$,} Proc.
London Math. Soc., II. Ser., 40 (1935), 336-344.

\item{[Wa]} T. Watson,
 {Rankin triple products and quantum chaos}, {\it Thesis, Princeton University, 2002}

\bye